\input amssym.def
\magnification=\magstep1
\baselineskip=13pt
\def \C {{\Bbb C}}
\def \Q {{\Bbb Q}}
\def \ZZ {{\Bbb Z}}
\def \Spec {{\rm Spec\,}}
\def\sdp{{\times\kern-.27em\raise1pt
  \hbox{$\scriptscriptstyle{|}$}
  {\kern-.25em\raise1pt\hbox{$\scriptscriptstyle{|}$}}
  \kern.1666em}}

\def \PP {{\Bbb P}}
\def \AA {{\Bbb A}}
\def \Gal {{\rm Gal}}
\def \Ind {{\rm Ind}}
\def \O {{\cal O}}

\def \H {{\rm H}}

\def \E {{\cal E}}
\def \ker {{\rm ker}}
\def \Ind {{\rm Ind}}
\def \qed {\hfill\lower0.9pt\vbox{\hrule \hbox{\vrule height 0.2 cm   
  \hskip 0.2 cm \vrule height 0.2 cm}\hrule}}
\def\cond#1{\par\noindent\rlap{#1}
  \hskip17pt\hangindent=20pt
            \hangafter=1}
\font\small=cmr8
\font\smit=cmti8
\font\smtt=cmtt8
\def \new {}

\centerline{\bf Local Galois theory in dimension two: Second edition} 

\centerline{David Harbater\footnote{$^*$}{\baselineskip=10pt
\small Supported in part by
NSF Grant DMS0200045.}
and Katherine F.~Stevenson\footnote{$^{**}$}
{\baselineskip=10pt
\small Supported in part by
NSA Grant ODOD-MDA 9049910038.
{\parindent=0pt\item{}  
{\smit 2000 Mathematics Subject Classification}.  Primary 12E30, 12F10, 14H30; Secondary 11S20, 12F12, 14J20.  
\item{} \baselineskip=10pt 
{\smit Key words and phrases}: absolute Galois group, embedding problem, fundamental group, Galois covers, patching, local, surfaces.}}}

\bigskip

{\narrower\noindent{\bf Abstract.} This paper proves a generalization of Shafarevich's Conjecture to fields of Laurent series in two variables over an arbitrary field.  This result says that the absolute Galois group $G_K$ of such a field $K$ is {\new semi-free} of rank equal to the cardinality of $K$, i.e.\ every non-trivial finite split embedding problem for $G_K$ has exactly ${\rm card}\,K$ proper solutions{\new, and moreover there is a set of this many {\it independent} solutions}.  {\new This strengthens the main result of the first edition of this paper, which appeared in 2005 and concerned quasi-free groups.}
We also strengthen a result of Pop and Haran-Jarden on the existence of proper regular solutions to split embedding problems for curves over large fields; our strengthening concerns integral models of curves, which are two-dimensional. \par}

\baselineskip=13pt

\medskip

\noindent  {\bf  Preface.}

\medskip

This manuscript is the second edition of a paper [HS05] that appeared in 2005, which introduced the notion of a quasi-free profinite group and studied this notion in the context of Galois theory.  Since then, various authors have used this notion to study absolute Galois groups.  In these applications, it has sometimes been useful to have a somewhat stronger notion, which at first we had called ``strongly quasi-free'' but which has now come to be known as ``semi-free''.  In those contexts, it would have been useful for authors to draw on versions of the results of [HS05] in the semi-free case.  
The proofs in [HS05] for quasi-free groups do indeed carry over to semi-free groups, though more care needs to be taken in various places and some arguments become more technical.  Due to the interest in this stronger notion, we produced a version of our manuscript that makes those changes, so that other authors could build upon the semi-free versions of the assertions.  We are now making this more available, in the form of a second edition to [HS05].  This edition does not attempt to bring up to date the discussions of related results (as we would in a new manuscript), except in a few places where it seemed especially important to do so.

\medskip

\noindent  {\bf  Section 1. Introduction.}

\medskip

This paper concerns the absolute Galois group of a field of the form $K = k((x,t))$, where $k$ is an arbitrary field.  Here $K$ is by definition the fraction field of the power series ring $k[[x,t]]$, the complete local ring of a point on a smooth surface.  
The absolute Galois group $G_K$ of $K$ can never be free, but we show that it is {\new semi-free} of rank equal to the cardinality of $K$.  This is quite different from the one-variable Laurent series case, where the absolute Galois groups are very far from being free.

Specifically, we prove the following result (see Theorem~5.1):

\medskip

\noindent{\bf Theorem 1.1.} {\sl If $k$ is a field, 
then the absolute Galois group of $K := k((x,t))$ is {\new semi-free} of rank ${\rm card}\, K$.}

\medskip

This theorem relies on a generalization of a result of Pop and Haran-Jarden on split embedding problems for curves over large fields (see Theorem~4.1 and the discussion below).  The notion of ``semi-free'', and a slightly weaker notion of ``quasi-free'', are introduced and discussed in Section~2.  There we show in particular (Theorem~2.1) that a profinite group of infinite rank is free if and only if it is projective and quasi-free; this is also equivalent to being projective and semi-free.

{\new If a group is semi-free (or quasi-free), then in particular it has the property} that every finite split embedding problem has a proper solution.  Thus Theorem~1.1 provides evidence for

\medskip

\noindent{\bf Conjecture 1.2.} (See [DD, \S2.1.2].) {\sl If $F$ is any Hilbertian field, then its absolute Galois group $G_F$ has the property that every finite split embedding problem has a proper solution.}

\medskip

Note that $K = k((x,t))$ is Hilbertian by a result of Weissauer [FJ, Theorem~14.17].)  Conjecture~1.2 is already known for Hilbertian fields that are large, by a result of Pop [Po96, Main Theorem~B].  {\new This raises the question of whether our field $K$ is large. In response to this question, there is a recent unpublished manuscript of Pop saying that it is.}  

For $K = k((x,t))$, the absolute Galois group $G_K$ cannot be free because it is not even projective (having cohomological dimension $>1$; cf.\ [AGV, Exp.~X, Cor.~2.4] and [Se, I, 3.4, Proposition~16]).  But by Theorem~2.1, our result says that $G_K$ is ``as close as possible to being free,'' given that it is not projective.  The lack of projectivity means that not every finite embedding problem has even a {\it weak} solution.  So in Theorem~1.1, the restriction to {\it split} embedding problems is essential.  Of course every split embedding problem has a weak solution, induced by the splitting; and our result asserts that there are {\it proper} solutions.  (See Section~2 for a review of definitions of these terms.)  

Theorem~1.1 can be regarded as a higher dimensional local version of Shafarevich's Conjecture (cf.~[Ha02]).  That conjecture says that the absolute Galois group of $\Q^{\rm ab}$ is a free profinite group of countable rank, and more generally that the absolute Galois group of the maximal cyclotomic extension of a global field is free of countable rank.  This assertion is known in the geometric case, where more generally it has been shown that the absolute Galois group of the function field of a curve over an algebraically closed field is free ([Ha95], [Po95]).  
This absolute Galois group was previously known to be projective; so in the 
terminology of the current paper, it sufficed to show that it was also quasi-free.  Thus we may regard a more general form of Shafarevich's Conjecture (even in higher dimensions, where projectivity fails) as saying that the absolute Galois groups of the function fields of certain schemes are quasi-free{\new; or in an even stronger form that it is semi-free}.  This remains open in the global case in dimension $>1$, but here we show it in the smooth equidimensional local case in dimension $2$.

Our result implies in particular that every finite group is a Galois group over $K = k((x,t))$.  This consequence was previously proven by T.~Lefcourt [Le], using that every finite group is the Galois group of a regular cover of the $K$-line and the fact that $K$ is Hilbertian.  Also, in the special case that $k=\C$, Theorem~5.3.9 of [Ha03] proved that every finite split embedding problem for $K$ has a proper solution (though the number of solutions was not considered there).  

One may wish to consider $K = k((x,t))$ as a ``two-dimensional local field'', but the situation differs from the one-variable case.  In dimension one, the notion of being local is essentially unambiguous, indicating the fraction field of a complete discrete valuation ring.  In dimension two, we can correspondingly consider the fraction field of a two-dimensional complete regular local ring, in particular the fraction field $K$ of $k[[x,t]]$.  But $K$ is also the fraction field of non-local rings, e.g.\ of 
$k[[x,t]][1/t]$, a Dedekind domain with infinitely many maximal ideals (corresponding to the height one primes of $k[[x,t]]$ other than $(t)$).  Similarly, $K$ is the function field of the blow-up of $\Spec k[[x,t]]$ at the closed point; and this blow-up is highly non-local, 
containing a copy of the projective line as
the exceptional divisor.  This aspect makes the higher dimensional situation different from the one-dimensional theory.  At the same time, there exists a ``more highly local'' two dimensional field, viz.\ the iterated Laurent series field $k((x))((t))$, which strictly contains $K = k((x,t))$.  But this is a discrete valuation field over $k((x))$, and its Galois theory is better understood (e.g.\ see [HP] in the case that $k$ is algebraically closed), with its absolute Galois group being very far from free.  In still higher dimensions, there is a whole ``zoo'' of fields between $k((x_1,\dots,x_n))$ and $k((x_1))\cdots((x_n))$ that are increasingly ``local'' and having increasingly simpler absolute Galois groups.  The case we consider here is thus the first of a much larger class.

In the process of proving Theorem~1.1, we also establish a strengthening of a result of Pop and Haran-Jarden.  The result of Pop (see [Po96, Main Theorem~A] and [Ha03, Theorem 5.1.9]) says that for any large field $F$, every finite split embedding problem for a one-variable function field over $F$ has a proper regular solution.  It is this result of Pop that implied Conjecture~1.2 in the case of large Hilbertian fields.  Haran and Jarden [HJ, Theorem~6.4] strengthened Pop's result to say that the proper regular solution may be chosen to be totally split over a given unramified point.  Here we show (Theorems~4.1 and 4.3) that there are infinitely many such solutions, which can be chosen to be totally split over any given finite set of (possibly ramified) points.  Moreover we show that if $F = k((t))$ and if we are given a smooth model for the curve over $k[[t]]$, then the solution may be chosen to be totally ramified over a given point on the closed fibre of the given cover, and there are exactly ${\rm card}\,F$ such solutions.  In doing so, it suffices to show that there are at least ${\rm card}\,F$ such solutions, since the other inequality follows from the fact that there are at most ${\rm card}\,F$ covers of the given $F$-curve.

The proof of Theorem~1.1, like that of the special case of $k=\C$ in [Ha03], uses formal patching and blowing up.  But while the proof in the case $k=\C$ could rely on the fact that $\C$ is algebraically closed of characteristic $0$, here the proof must proceed differently, since $k$ can now be an arbitrary field.  Instead, we use an argument that relies on the strengthened version of the theorem of Pop and Haran-Jarden described above.
  
\medskip

The structure of this paper is as follows:  In Section~2, which contains foundational material about embedding problems, we define and study the notion of a ``quasi-free'' {\new and ``semi-free''} profinite group of infinite rank.  In particular, we show that being free is equivalent to being quasi-free and projective; and we also prove results in the countably generated case (Proposition~2.7, Corollary~2.8) which are related to Iwasawa's theorem [Iw, p.567].  Section~4 contains the strengthened version of the result of Pop and Haran-Jarden, with Section~3 containing preliminary results for that proof.  
Finally, in Section~5 we prove Theorem~1.1 above (Theorem~5.1 there), that $G_K$ is {\new semi-free} of rank equal to ${\rm card}\,K$.

\medskip

The techniques in this paper include the deformation of degenerate covers --- an approach that was used in previous work of the present authors, including their Ph.D.\ theses (written under the supervision of M.~Artin and the first author, respectively).  The authors wish to acknowledge Michael Artin's inspiration for this circle of ideas.

\bigskip

{\it Terminology:}  A (branched) {\it cover} of schemes $f:Y \to X$ is a finite generically separable morphism.  If $f:Y \to X$ is a cover, then the {\it Galois group} $\Gal(Y/X)$ consists of the automorphisms $\phi$ of $Y$ satisfying $f\phi = f$.  If $G$ is a finite group, then a $G$-{\it Galois cover} is a cover $Y \to X$ together with an injection $G \hookrightarrow \Gal(Y/X)$ such that $\O_Y^G = f^*(\O_X)$ (where the left side denotes the sheaf of $G$-invariants).  If $X$ is a normal scheme, this condition is equivalent to saying that $G$ acts simply transitively on a generic geometric fibre of $Y \to X$.  A cover $Y \to X$ is {\it connected} [resp.\ {\it irreducible, normal}], if the schemes $X$ and $Y$ are.

Let $f:Y \to X$ be a $G$-Galois cover, let $Q$ be a point of $Y$ and let $P = f(Q) \in X$.  The {\it decomposition group} at $Q$ is the subgroup of $G$ consisting of elements that fix the point $Q$.  The {\it inertia group} at $Q$ is the subset of the decomposition group that induces the identity automorphism on the residue field at $Q$.  The cover is {\it ramified} at $Q$ if the inertia group is non-trivial, and it is {\it totally ramified} at $Q$ if the inertia group is all of $G$ (in which case $Q$ is the only point in its fibre over $P$).  The cover is {\it totally split} (or {\it splits completely}) over a point $P \in X$ if each decomposition group over $P$ is trivial.  

Let $k$ be an arbitrary field and let $X$ be an integral $k$-scheme such that $k$ is algebraically closed in its function field.  Let $Y \to X$ be a connected cover.  We say that $Y \to X$ is a {\it regular} cover of $k$-schemes if $Y$ is geometrically connected as a $k$-scheme (or equivalently, if $k$ is algebraically closed in the function field of $Y$).   We say that $Y \to X$ is a {\it purely arithmetic} cover of $k$-schemes if $Y$ is isomorphic to $X \times_k \ell$ as a cover of $X$, where $\ell$ is a finite (necessarily separable) field extension of $k$.  Thus for any cover $Y \to X$, if $\ell$ is the algebraic closure of $k$ in the function field of $Y$, then $Y \to X$ factors as $Y \to X_\ell \to X$, where $X_\ell = X \times_k \ell$; here $Y \to X_\ell$ is regular and $X_\ell \to X$ is purely arithmetic.  (Note that this notion of ``regular cover'' is unrelated to the notion of a ``regular scheme'' in the sense of having regular local rings.)

If $H \subset G$ and if $Y \to X$ is an $H$-Galois cover, then the {\it induced} $G$-Galois cover $\Ind_H^G Y \to X$ is the $G$-Galois cover of $X$ obtained by taking a disjoint union of copies of $Y \to X$, indexed by the left cosets of $H$ in $G$.  Here the stabilizer of
the ``identity copy'' of $Y$ (corresponding to the identity coset) is $H$, and the stabilizers of the other copies are the conjugates of $H$ in $G$.  The {\it trivial} $G$-Galois cover of $X$ is $\Ind_1^G X \to X$; this consists of disjoint copies of $X$ that are indexed by the elements of $G$ and are permuted according to the regular representation.

If $X$ is irreducible, then a cover $f:Y \to X$ is a {\it mock cover} if the restriction of $f$ to each irreducible component of $Y$ is an isomorphism to $X$.  

If $p$ is a prime number, then a {\it cyclic-by-$p$} group is a semi-direct product $P \sdp C$, where $P$ is a finite $p$-group and $C$ is a cyclic group of order prime to $p$.  If $p=0$, then by a $p$-{\it group} we will mean the trivial group, and by a {\it cyclic-by-$p$} group we will mean a cyclic group.  Thus if $k$ is an algebraically closed field of characteristic $p \ge 0$, and if $\hat \ell$ is a finite Galois field extension of $\hat k = k((t))$, then the associated Galois group $\Gal(\hat \ell/\hat k)$ is cyclic-by-$p$.

\bigskip

\noindent {\bf Section 2. Embedding Problems and Quasi-Free Profinite Groups}

\medskip

In this section we define the notions of ``quasi-free'' {\new and ``semi-free''} profinite groups, and prove that a profinite group of infinite rank is free if and only if it is quasi-free and projective.  {\new These notions are} defined in terms of embedding problems for profinite groups.  We also prove related results in the countably generated case (Proposition~2.7 and Corollary~2.8).  In addition, we
review basic definitions, discuss the relationship of embedding problems to Galois theory, and prove a technical result (Proposition~2.9) which will be useful in Section~4.  (See also Chapters 24 and 25 of the 2005 second edition of [FJ] for a further discussion along these lines.)

If $\Pi$ is a profinite group and $m$ is an infinite cardinal number, then we will say that $\Pi$ is {\it quasi-free} of {\it rank} $m$ if every non-trivial finite split embedding problem for $\Pi$ has exactly $m$ proper solutions.  {\new We also say that $\Pi$ is {\it semi-free} of rank $m$ if} $\Pi$ has rank $m$ and every non-trivial finite split embedding problem $\E$ for $\Pi$ has a set of $m$ independent proper solutions.  (See below for definitions concerning profinite groups and embedding problems.)  {\new These conditions generalize} the condition of being free of rank $m$, to the class of profinite groups that are not necessarily projective.  We will prove:

\medskip

\noindent{\bf Theorem~2.1.} {\sl Let $\Pi$ be a profinite group and let $m$ be an infinite cardinal.  Then $\Pi$ is a free profinite group of rank $m$ if and only if the following conditions are satisfied:

(i) $\Pi$ is projective.

(ii) $\Pi$ is quasi-free of rank $m$.}

\medskip

{\new Notice that the analog of Theorem 2.1 for semi-free profinite groups follows immediately from Theorem 2.1 itself as semi-free implies quasi-free, and both are implied by free.}

\medskip

\noindent{\bf Remark~2.2.}  This theorem is a variant on a result of Melnikov and Chatzidakis [Ja, Lemma~2.1].  That result says that a profinite group $\Pi$ is free of rank $m$ if and only if every (not necessarily split) non-trivial finite embedding problem for $\Pi$ has exactly $m$ proper solutions.  Like the theorem above, the key direction is the {\new reverse} implication.  Melnikov had proved that direction of the earlier result under the additional assumption that $\Pi$ has rank $m$ (see [FJ, Prop.~24.18]); and Chatzidakis later showed that in fact the rank is automatically $m$.

\medskip

Before proving the theorem, we recall some basic definitions and facts about profinite groups and embedding problems (cf.\ [FJ]), and prove some preliminary results.  

In the category of profinite groups, all {\it homomorphisms} are required to be continuous, and {\it generating sets} are taken in the profinite sense (i.e.\ the generating condition is that there are no proper closed subgroups containing the set).  A generating set $S$ of a profinite group $\Pi$ {\it converges to $1$} if every open normal subgroup of $\Pi$ contains all but finitely many elements of $S$.  Such a generating set always exists, by a result of Douady [FJ, Prop.~15.11].  The minimal cardinality of such a generating set is called the {\it rank} of $\Pi$.  In the case that $\Pi$ is finitely generated (as a profinite group), the rank is the minimal number of generators of $\Pi$.  If $\Pi$ is not finitely generated, then the rank is the cardinality of any generating set $S$ that converges to $1$; this is independent of the choice of $S$, by [FJ, Supplement~15.12].  In fact that result says that ${\rm card}\,S$ (and hence the rank of $\Pi$) is equal to the cardinality of the set of all open normal subgroups of $\Pi$.  

A profinite group $\Pi$ is {\it free} if there is a generating set $S$ that converges to $1$ and that has the following additional property: for every profinite group $\Delta$ and every map $\phi_0:S \to \Delta$ such that $\phi_0(S) \subset \Delta$ converges to $1$, there exists an extension of $\phi_0$ to a homomorphism $\phi:\Pi \to \Delta$.  For every cardinal $m$, there is (up to isomorphism) a unique free profinite group of rank $m$, denoted by $\hat F_m$ [FJ, p.191].  

An {\it embedding problem} $\E$ is a pair $(\alpha:\Pi \to G, f:\Gamma \to G)$ of epimorphisms of profinite groups.  The {\it kernel} of $\E$ is $\ker(f)$.  We say that  $\E$ is {\it finite} if $\Gamma$ is finite; it is {\it non-trivial} if its kernel is non-trivial; and it is split if $f$ has a section.  A {\it weak solution} to $\E$ is a homomorphism $\lambda:\Pi \to \Gamma$ such that $f\lambda=\alpha$.  A {\it proper solution} to $\E$ is a weak solution in which $\lambda$ is surjective.  A profinite group $\Pi$ is {\it projective} if every finite embedding problem has a weak solution.  (Note that for any profinite group $\Pi$, every finite split embedding problem has a weak solution.)

In fact, a profinite group $\Pi$ is projective if and only if it is isomorphic to a closed subgroup of a free profinite group [FJ, Cor.~20.14]; in particular, every free profinite group is projective.  Also, a profinite group $\Pi$ is projective if and only if it has cohomological dimension $\le 1$, by [Gru, Theorem~4] (or by [Se, I, \S3.4 Prop.~16 and \S5.9 Prop.~45], using that ${\rm cd} = {\rm max}\ {\rm cd}_p$).

{\new Given a finite split embedding problem $\E$ for a profinite group $\Pi$,} {\new we say that a set of solutions is {\it independent} if for every finite subset $\lambda_1,\dots,\lambda_n$, their kernels $N_i := \ker \lambda_i$ satisfy the condition $(F:\cap_{i=1}^n N_i) = \prod_{i=1}^n (F:N_i)$.  If $\E = (\alpha:G_K \to G, f:\Gamma \to G)$ is an embedding problem for a field $K$ (in the sense given in the discussion before Proposition~2.9 {\new below}), then a set of proper solutions is independent if and only if the corresponding field extensions are linearly disjoint over $L$, the $G$-Galois field extension of $K$ corresponding to $\alpha$.  By the theorem of Melnikov and Chatzidakis referred to in Remark~2.2, {\new semi-freeness} is equivalent to the condition that every non-trivial finite split embedding problem $\E$ for $\Pi$ has exactly $m$ proper solutions, and that these solutions include a set of $m$ independent proper solutions.}  

\medskip

For any two profinite groups $\Pi, \Delta$, let ${\rm Epi}(\Pi,\Delta)$ be the set of epimorphisms $\Pi \to \Delta$.

\medskip

\noindent{\bf Lemma~2.3.} {\sl Let $\Pi$ be a profinite group.  Suppose that ${\rm Epi}(\Pi,G)$ is non-empty for every finite group $G$.  Then the rank of $\Pi$ is the sum of the cardinalities of the sets ${\rm Epi}(\Pi,G)$, where $G$ ranges over isomorphism classes of finite groups.}

\medskip

\noindent{\it Proof.}  First observe that $\Pi$ is not a finitely generated profinite group.  Namely, if it were generated by a set $S$ of $n$ elements, then ${\rm Epi}(\Pi,G)$ would be empty if $G$ is the product of $n+1$ copies of the cyclic group of two elements.  

So [FJ, Supplement~15.12] applies, and says that the rank of $\Pi$ is equal to the cardinality of the set of open normal subgroups of $\Pi$.  But this is the same as the cardinality of the set of epimorphisms from $\Pi$ to finite groups $G$ (up to isomorphism), i.e.\ the sum of the cardinalities of ${\rm Epi}(\Pi,G)$, where $G$ ranges over isomorphism classes of finite groups.   \qed

\medskip

Let $\Pi$ be a profinite group and let $\E = (\alpha:\Pi \to G, f:\Gamma \to G)$ be a finite embedding problem for $\Pi$.  Suppose that the epimorphism $\alpha:\Pi \to G$ factors as $r\alpha'$, where $\alpha':\Pi \to G'$ and $r:G'\to G$ are 
epimorphisms, for some finite group $G'$.  We consider the {\it induced embedding problem} $\E_{\alpha'} = (\alpha':\Pi \to G', f':\Gamma' \to G')$ by taking $\Gamma' = \Gamma \times_G G'$ and letting $f':\Gamma' \to G'$ be the second projection map.  Here 
$f'$ is surjective because $f$ is; and so $\E_{\alpha'}$ is a finite embedding problem.  Here $\E$ and $\E_{\alpha'}$ have isomorphic kernels; indeed $\ker(f') = \ker(f) \times 1 \subset \Gamma \times_G G' = \Gamma'$.  Note also that the first projection map $q:\Gamma' = \Gamma \times_G G' \to \Gamma$ is surjective since $r:G' \to G$ is surjective; and $f q = r f'$.   

In this situation, every proper solution $\lambda':\Pi\to \Gamma'$ of $\E_{\alpha'}$ induces a proper solution $\lambda := q\lambda':\Pi \to \Gamma$ of $\E$; viz.\ $f\lambda = f q \lambda' = r f' \lambda' = r\alpha' = \alpha$, and 
$\lambda$ is surjective because $q$ and $\lambda'$ are.  So we obtain a map ${\rm PS}(\E_{\alpha'}) \to {\rm PS}(\E)$, where $\rm PS$ denotes the set of proper solutions to the embedding problem.  

Consider a non-trivial finite split embedding problem $\E = (\alpha:\Pi \to G,\, f:\Gamma \to G)$ for $\Pi$.
For any positive integer $n$, let $\Gamma^n_G$ denote the $n^{\rm th}$ fibre power of $\Gamma$ over $G$, and let $f^n: \Gamma^n_G \to G$ be the natural projection map (i.e.\ the composition of $f$ with the $i^{\rm th}$ projection map $\Gamma^n_G \to G$, for any $i$).  Note that $\Gamma'^n_{G'}$ is naturally isomorphic to $\Gamma^n_G \times_G G'$.
{\new Observe that} $\E^n = (\alpha:\Pi \to G,\, f^n:\Gamma^n_G \to G)$ is also a non-trivial finite embedding problem, having a splitting given by composing the splitting of $\E$ with the diagonal embedding $\Gamma \hookrightarrow \Gamma^n_G$.  If $\beta_1,\dots,\beta_n:\Pi \to \Gamma$, then
$\beta:= (\beta_1,\dots,\beta_n)$ defines a proper regular solution to $\E^n$  {\new if and only if $\beta_1,\dots,\beta_n$ are independent proper regular solutions to $\E$}.  

\medskip

\noindent{\bf Lemma~2.4.} {\sl In the above situation, the map ${\rm PS}(\E_{\alpha'}) \to {\rm PS}(\E)$ is injective.  Moreover, under this map, a set of independent proper solutions is mapped to a set of independent proper solutions.}

\medskip

\noindent{\it Proof.}  We have $\E = (\alpha:\Pi \to G, f:\Gamma \to G)$ and $\E_{\alpha'} = (\alpha':\Pi \to G', f':\Gamma' \to G')$.
Say $\lambda_1', \lambda_2':\Pi \to \Gamma'$ are proper solutions to $\E_{\alpha'}$.  So $f'\lambda_1' = \alpha' = f'\lambda_2'$.   If $\lambda_1'$ and $\lambda_2'$ have the same image under ${\rm PS}(\E_{\alpha'}) \to {\rm PS}(\E)$, then $q\lambda_1' = q\lambda_2'$.  So $\lambda_1', \lambda_2'$ have the same composition with $(q,f'):\Gamma' \to \Gamma \times_G G' = \Gamma'$.  But $(q,f')$ is the identity map on $\Gamma'$.  So $\lambda_1' = \lambda_2'$.  This proves injectivity.

{\new To complete the proof, it remains to show that a set of proper solutions to $\E$, induced by a set of {\new independent} proper solutions to $\E_{\alpha'}$, is itself {\new independent}.  For any positive integer $n$, 
{\new consider the $n$th fibre powers $\Gamma^n_G$ and $\Gamma'^n_{G'}$, and the associated finite split embedding problems $\E^n$ and $\E^n_{\alpha'}$} as discussed above.  
By the first part of the proof, the map ${\rm PS}(\E^n_{\alpha'}) \to {\rm PS}(\E^n)$ is injective.  Thus if $\beta_i$ is the image of $\beta_i'$ under ${\rm PS}(\E_{\alpha'}) \to {\rm PS}(\E)$, then $\beta := (\beta_1,\dots,\beta_n)$ defines a proper regular solution to $\E^n$.  Hence $\beta_1,\dots,\beta_n$ are {\new independent} proper solutions to $\E$, as desired.}
\qed

 \medskip

\noindent{\bf Example~2.5.} In certain key cases the above map $r:G' \to G$ factors as $f\phi$ for some $\phi:G' \to \Gamma$:

(a) If $G' \subset \Gamma$ and $r = f|_{G'}$, then $\alpha'$ is just a weak solution to $\E$.  (Conversely, any weak solution $\alpha'$ to $\E$ induces such a $G' := {\rm image}(\alpha')$ and an $r = f|_{G'}$ with $\alpha=r\alpha'$.)  In this case we may take $\phi$ to be the inclusion $G' \hookrightarrow \Gamma$.

(b) In the general situation considered above, if $\E$ is a split embedding problem and if $s:G \to \Gamma$ is a section for $f$, then we may take $\phi=sr$.

\medskip

\noindent{\bf Lemma 2.6.} {\sl In the general situation above, suppose that $r:G' \to G$ factors through $f:\Gamma \to G$, say as $r = f\phi$ with $\phi:G' \to \Gamma$.  Then the induced embedding problem $\E_{\alpha'}$ has a splitting $s':G' \to \Gamma'$, given by $s' = (\phi,{\rm id}_{G'})$.
In particular, $\E_{\alpha'}$ is split in the situations of Examples~2.5 (a) and (b) above.} 

\medskip

\noindent {\it Proof.} With $s'$ as above, $f's' = {\rm id}_{G'}$, so $s'$ is a splitting.  \qed

\medskip 

Note that in the situation of Example~2.5(a), we may identify the splitting $s'$ with the diagonal map $G' \to G' \times_G G' \subset \Gamma \times_G G' = \Gamma'$.  In the situation of Example~2.5(b), $s'$ lifts $s$ in the sense that $qs' = sr$, because $qs' = q \circ (\phi,{\rm id}_{G'}) = q \circ (sr,{\rm id}_{G'}) = sr$.

\medskip

We can now prove Theorem~2.1 above:

\medskip

\noindent{\it Proof of Theorem~2.1.} The forward direction follows from [FJ, Cor.~20.14] (for (i)) and [FJ, Lemma~24.14] (for (ii)).

For the reverse direction, let $G$ be a non-trivial finite group, and consider the non-trivial finite split embedding problem $\E_G := (\Pi \to 1, G \to 1)$.  By (ii), this has exactly $m$ proper solutions.  That is, there are exactly $m$ epimorphisms $\Pi \to G$.  So by Lemma~2.3 above, the profinite group $\Pi$ has rank $m$ (using that $m$ is infinite and that there are countably many isomorphism classes of finite groups $G$).  It remains to show that $\Pi$ is {\it free} of rank $m$.

Let $\E$ be any non-trivial finite embedding problem for $\Pi$.  Since $\Pi$ is projective by (i), there is a weak solution $\alpha'$ to $\E$, and hence an induced non-trivial finite split embedding problem $\E_{\alpha'}$ for $\Pi$, as in Example~2.5(a).  By (ii), $\E_{\alpha'}$ has $m$ proper solutions.  So by Lemma~2.4, $\E$ has at least $m$ proper solutions.  But by Lemma~2.3, $\E$ has at most $m$ proper solutions, since $\Pi$ has rank $m$, and since every proper solution is an epimorphism from $\Pi$ to a fixed finite group.  So $\E$ has exactly $m$ proper solutions.  

We now conclude the proof using Melnikov's theorem [FJ, Prop.~24.18].  Namely, that result says that two profinite groups of the same infinite rank $m$ must be isomorphic provided that each has the property that every non-trivial finite embedding problem has exactly $m$ solutions.  As just shown, the profinite group
$\Pi$ has this property.  But so does the free profinite group of rank $m$, by [FJ, Lemma~24.14].  So $\Pi$ is isomorphic to this free profinite group.  \qed

\medskip

In the case that $m$ is countable, the quasi-free condition becomes a bit simpler:

\medskip

\noindent{\bf Proposition~2.7.} {\sl {\new A profinite group $\Pi$ of countably infinite rank is quasi-free} if and only if every finite split embedding problem for $\Pi$ has a proper solution.}

\medskip

\noindent{\it Proof.}  The forward direction is trivial.  For the reverse direction, consider a non-trivial finite split embedding problem $\E = (\alpha:\Pi \to G,\, f:\Gamma \to G)$ for $\Pi$.
For any positive integer $n$, {\new consider the finite split embedding problem 
$\E^n = (\alpha:\Pi \to G,\, f^n:\Gamma^n_G \to G)$, as in the paragraph before Lemma~2.4.}  By hypothesis, $\E^n$ has a proper solution $\lambda_n$.  Composing $\lambda_n$ with the $n$ projection maps $\Gamma^n_G \to \Gamma$ (in turn) yields $n$ distinct (in fact independent) proper solutions to $\E$.  Since this holds for every $n$, it follows that the set of proper solutions to $\E$ is infinite.  
But $\Pi$ has countable rank.  So the set of epimorphisms $\Pi \to \Gamma$ is at most countable.  So in fact the set of proper solutions to $\E$ is countably infinite, by Lemma~2.3 above.  Since this is true for all $\E$, the profinite group $\Pi$ is quasi-free of countably infinite rank. \qed

\medskip

As a consequence, we obtain the following result, which is related to Iwasawa's theorem ([Iw, p.567]; see below): 

\medskip

\noindent{\bf Corollary~2.8.} {\sl Let $\Pi$ be a profinite group of countably infinite rank.  Then $\Pi$ is a free profinite group (of countable rank) if and only if the following conditions are satisfied:

(i) $\Pi$ is projective.

(ii) Every finite split embedding problem for $\Pi$ has a proper solution.}

\medskip

\noindent{\it Proof.}  By Proposition~2.7, condition (ii) in the corollary is equivalent to condition (ii) of Theorem~2.1, in this situation.  So the corollary follows from Theorem~2.1.  \qed

\medskip

Note that the argument proving Proposition~2.7 remains valid if one instead considers {\it all} finite embedding problems, rather than just split ones.  That is, if a profinite group $\Pi$ has the property that every finite embedding problem for $\Pi$ has a proper solution, then every non-trivial finite embedding problem for $\Pi$ must have infinitely many proper solutions.  And so if $\Pi$ has countable rank, then each such embedding problem has exactly countably infinite proper solutions. 

This shows that the result of Melnikov and Chatzidakis generalizes the earlier related result of Iwasawa ([Iw, p.567]; cf.\ also [FJ, Cor.~24.2]) in the countable rank case.  Iwasawa's result says that a profinite group of countable rank is free if and only if every (not necessarily split) finite embedding problem has a proper solution.  

On the other hand, by relying on Iwasawa's result instead of using Theorem~2.1, we obtain another proof of the above corollary:  

\medskip

\noindent{\it Alternative proof of Corollary~2.8.}  The forward direction follows (as in the proof of Theorem~2.1) from [FJ, Cor.\ 20.14 and Lemma~24.14].  For the reverse direction, by Iwasawa's result we are reduced to proving that every finite embedding problem $\E$ for $\Pi$ has a proper solution.  For this, first note that $\E$ has a weak solution $\alpha'$ because $\Pi$ is projective.  The induced finite split embedding problem $\E_{\alpha'}$ has a proper solution $\lambda'$ by (ii).  The image of $\lambda'$ under ${\rm PS}(\E_{\alpha'}) \to {\rm PS}(\E)$ is then a proper solution to $\E$.  \qed

\medskip

Embedding problems arise in Galois theory by taking the profinite group $\Pi$ to be the absolute Galois group $G_K$ of some field $K$.  If $K$ is a field, then by a {\it finite embedding problem} $\E = (\alpha:G_K \to G, f:\Gamma \to G)$ {\it for} $K$ we will mean such an embedding problem for $G_K$.  Giving such a problem corresponds to giving a $G$-Galois field extension $L$ of $K$, where $G$ is a quotient of $\Gamma$; and a proper solution to $\E$ corresponds to giving a $\Gamma$-Galois field extension $M$ of $K$ that contains $L$.  That is, giving a proper solution is equivalent to embedding the given $G$-Galois field extension of $K$ into a $\Gamma$-Galois field extension of $K$, compatibly with the quotient map $f:\Gamma \to G$ (and this is the origin of the terminology ``embedding problem'').  
If $k$ is a subfield of $K$ and we regard $K$ as a $k$-algebra, then we say that a proper solution as above is {\it regular} (over $k$) if the algebraic closures of $k$ in $L$ and in $M$ are the same (viewing $L \subset M$).  

In particular, given a field $k$ and a connected normal $k$-scheme $X$ with function field $K$, a finite embedding problem $\E = (\alpha:G_K \to G, f:\Gamma \to G)$ for $K$ corresponds to giving a $G$-Galois normal connected (branched) cover $Y \to X$.  A weak solution to $\E$ corresponds to giving a $\Gamma$-Galois normal cover $Z \to X$ that dominates $Y \to X$; such a solution is proper if and only if $Z$ is also connected.  A {\it proper regular solution} to $\E$ (over $k$) is a proper solution such that the algebraic closures of $k$ in the function fields of $Y$ and $Z$ are the same; here we regard the function field of $Y$ as contained in that of $Z$.  (This notion generalizes the notion of a ``regular cover'' of $k$-schemes, by considering the case in which $G$ is trivial.)

\medskip

With notation as before we have the following result:

\medskip

\noindent{\bf Proposition~2.9}.  {\sl Let $K$ be a field and let $\E = (\alpha:G_K \to G, f:\Gamma \to G)$ be a finite embedding problem for $K$.  Let $L$ be the $G$-Galois field extension of $K$ corresponding to $\alpha$.  Let $L'$ be a finite Galois extension of $K$ that contains $L$, say with Galois group $G'$ over $K$, and let $\alpha':G_K \to G'$ be the corresponding epimorphism.  Consider 
the induced embedding problem $\E_{\alpha'} = (\alpha':G_K \to G', f':\Gamma' \to G')$, where $\Gamma' = \Gamma \times_G G'$.  Let $\lambda':G_K \to \Gamma'$ be a proper solution to $\E_{\alpha'}$, let $\lambda:G_K \to \Gamma$ be the image of $\lambda'$ under ${\rm PS}(\E_{\alpha'}) \to {\rm PS}(\E)$, and let $M'$ and $M$ be the corresponding Galois field extensions of $K$, with groups $\Gamma'$ and $\Gamma$ respectively.   Then

a) $M'=ML'$.

b) $M \cap L'  = L$ in $M'$, and the natural map $M \otimes_L L' \to M'$ is an isomorphism.

c) If $\lambda'$ is a regular solution to $\E_{\alpha'}$ over $k$,§ then $\lambda$ is a regular solution to $\E$ over $k$.} 

\medskip

\noindent{\it Proof.} (a) Since $\Gamma' = \Gamma \times_G G'$, we have that $G_{M'} = \ker(\lambda') = \ker(\lambda) \cap \ker(
\alpha') = G_M \cap G_{L'}$.  So $M'=ML'$.

(b) The natural map $M \otimes_L L' \to M'$ is surjective by (a).  But $[M:L] = [M':L']$ since the kernels of $\E$ and $\E_{\alpha'}$ are isomorphic.  So $[M \otimes_L L':L] = [M:L] [L':L] = [M':L]$ and hence the map is an isomorphism.  Thus $M \cap L'  = L$ in $M'$, since this holds in the tensor product.

(c) By hypothesis, the algebraic closures of $k$ in $L'$ and in $M'$ are equal.  We wish to show that the algebraic closures of $k$ in $L$ and in $M$ are equal.  So suppose that $u \in M$ is algebraic over $k$.  Then $u \in M \subset M'$; so by hypothesis we also have $u \in L'$.  But by part (b), $L' \cap M = L$.  So actually $u \in L$.  This 
shows that the algebraic closure of $k$ in $L$ contains the algebraic closure of $k$ in $M$; and the other containment is trivial.
\qed

\bigskip

\noindent {\bf Section~3. Models of covers and blowing up}

\medskip

This section contains several results about models of covers of curves over complete discrete valuation rings $R$ (especially $R=k[[t]]$), for use in Section~4.  In particular we consider the effect on covers of blowing up at points on the closed fibre.  This will be useful later, in constructing $R$-models of covers with reducible closed fibres.

The reader may wish to skip this section initially, and to refer back to it as needed in Section~4.

\medskip

\noindent{\bf Proposition~3.1.} {\sl Let $R$ be a complete d.v.r. 
with fraction field $F$.  Let $X$ be a smooth projective connected curve over $F$.

a) Then there is a normal proper model $\bar X$ for $X$ over $R$ together with a finite $R$-morphism $\phi: \bar X \to \PP^1_R$ that is unramified at the generic point of $\bar X$.

b) Suppose that $\bar X$ is a smooth proper model for $X$ over $R$.  Then there is a finite $R$-morphism $\phi:\bar X \to \PP^1_R$ that is unramified at the generic point of the special fibre.}

\medskip

{\it Proof.}  Let $K$ be the function field of $X$ and let $k$ be the residue field of $R$.  

(a) Since $X$ is smooth over $F$, it follows that $K$ is separably generated over $F$; let $\{f\}$ be a separating transcendence basis for $K$ over $F$.  This yields a branched covering morphism $\phi_F:X \to \PP^1_F$.  Here we   
view $\PP^1_F$ as the generic fibre of $\PP^1_R$, and we let $\bar X$ be the normalization of $\PP^1_R$ in $X$.  So $\phi_F$ extends to a morphism $\phi:\bar X \to \PP^1_R$.
Since $K$ is finite and separable over $F(f)$, the morphism $\phi$ has the desired property.

(b) The closed fibre $X_0 = \bar X \times_R k$ of $\bar X$ is connected, by Zariski's Connectedness Theorem [Hrt, III, Cor.~11.3].  Since $X_0$ is 
smooth over $k$, it is irreducible, with one generic point $\eta$.  Note that a uniformizer for $R$ is also a uniformizer for the local ring $\O_{\bar X, \eta}$, since $X_0$ is reduced.  Also by smoothness, the function field $K_0$ of $X_0$ has a separating transcendance basis $\{f_0\}$ over $k$.   
 
Let $C_0, D_0$ be the zero and pole loci of the rational function $f_0$ on the smooth proper $k$-curve $X_0$.  Write $D_0 = \sum_{i=1}^r a_i P_{i,0}$ with $a_i > 0$, where the $P_{i,0}$'s are distinct closed points on $X_0$, say with local uniformizers $\pi_{i,0} \in \O_{X_0,P_{i,0}}$ on $X_0$.  Lifting $\pi_{i,0}$ to an element $\pi_i \in \O_{\bar X,P_{i,0}}$, we obtain an effective divisor $\bar P_i$ on $\bar X$ whose restriction to $X_0$ is $P_{i,0}$.  So $\bar D := \sum_{i=1}^r a_i \bar P_i$ is an effective divisor on $\bar X$ whose restriction to $X_0$ is $D_0$.  

We claim that the canonical homomorphism $\H^0(\bar X,\O(n\bar D)) \to \H^0(X_0,\O(nD_0))$ is surjective for all sufficiently large $n$.  To see this, let $\hat X$ be the formal scheme associated to $\bar X$, let $\hat D$ be the divisor on $\hat X$ associated to $\bar D$, and let $\O_0$ be the reduction of $\O := \O_{\hat X}$ modulo the maximal ideal of $R$.  If $n >\!\!> 0$, then the canonical homomorphism $\H^0(\hat X, \O(n\hat D)) \to \H^0(\hat X, \O_0(n\hat D))$ is a surjection [Gr61, Prop.~5.2.3].  The claim then follows by applying the canonical isomorphisms $\H^0(\hat X, \O(n\hat D)) = \H^0(\bar X,\O(n\bar D))$ [Gr61, Prop.~5.1.2] and
$\H^0(\hat X, \O_0(n\hat D)) = \H^0(X_0,\O(nD_0))$.

Choose such a sufficiently large $n$ that is not divisible by the characteristic of $k$.  So the rational function $f_0^n$ on the smooth $k$-curve $X_0$ still gives a separating transcendence basis for $K_0$ over $k$.  Thus the corresponding morphism $X_0 \to \PP^1_k$ is finite and generically separable.

The divisor of $f_0^n$ on $X_0$ is $nC_0 - nD_0$; so there is some $g \in \H^0(\bar X,\O(n\bar D))$ that maps to $f_0^n \in \H^0(X_0,\O(nD_0))$ under the above surjective map.  Let $(g)_0$, $(g)_\infty$ be the divisors of zero and poles of $g$ on $\bar X$.  Thus $(g)_0 - (g)_\infty$ restricts to $nC_0 - nD_0$, the divisor of $f_0^n$, on $X_0$.  Since $C_0, D_0$ have disjoint support, it follows that the restriction of $(g)_0$ [resp.\ $(g)_\infty$] to $X_0$ is at least $nC_0$ [resp.\ $nD_0$].  But by definition of $g$, its divisor of poles $(g)_\infty$ on $\bar X$ is at most $n\bar D$, whose restriction to $X_0$ is $nD_0$.  So in fact $(g)_\infty = n\bar D$, with restriction $nD_0$.  Thus the restriction of $(g)_0$ to $X_0$ is exactly $nC_0$.  Hence the supports of $(g)_0, (g)_\infty$ are disjoint on $\bar X$; and so the rational function $g$ on the smooth $R$-curve $\bar X$ has no locus of indeterminacy.  Thus $g$ defines a morphism $\phi:\bar X \to \PP^1_R$ over $R$.  Its restriction to the closed fibre is the morphism given by $f_0^n$, which is generically separable and hence unramified at the generic point $\eta$ of $X_0$ (and so also at the generic point of $\bar X$).  

The $R$-morphism $\phi$ is finite-to-one on the closed fibre since 
the fibre is irreducible and since the morphism
is unramified at the generic point of the fibre; and similarly it is finite-to-one on the generic fibre.  It is proper since $\bar X$ is proper over $R$.   Hence $\phi$ is finite.  So the morphism is as asserted.  \qed

\medskip

\noindent{\it Remark.} A related result appears at [GMP, Theorem~3.1].

\medskip

Let $R$ be a complete d.v.r.\ and let $\bar Y$ be an irreducible normal $R$-curve.  Let $Y$ be the general fibre of $\bar Y$, let $B$ be a reduced proper closed subset of $Y$, and let $\Sigma$ be the specialization of $B$ to the closed fibre of $\bar Y$ (i.e.\ the reduced intersection of the closure of $B$ in $\bar Y$ with the closed fibre of $\bar Y$).  Write $\bar Y_0 = \bar Y$ and $\Sigma_0 = \Sigma$, and let $\bar Y_1$ be the normalization of the blow-up of $\bar Y_0$ at the points of $\Sigma_0$.  So there is a birational morphism $\bar Y_1 \to \bar Y_0$ which induces an isomorphism between their general fibres (and so identifies $Y$ with the general fibre of $\bar Y_1$).  Regarding $B \subset Y \subset \bar Y_1$, let
$\Sigma_1$ be the specialization of $B$ to the closed fibre of $\bar Y_1$.  We call $(\bar Y_1, B, \Sigma_1)$ the {\it blow-up} of $(\bar Y, B, \Sigma)$.  Inductively, define the {\it $n^{\rm th}$ blow-up} $(\bar Y_n, B, \Sigma_n)$ of $(\bar Y, B, \Sigma)$ to be the blow-up of $(\bar Y_{n-1}, B, \Sigma_{n-1})$.  So $\bar Y_n$ is an irreducible normal $R$-curve that is birational to $\bar Y$ and whose general fibre is equipped with an isomorphism to $Y$.  For short, we will say that $\bar Y_n$ is the $n^{\rm th}$ blow-up of $\bar Y$ with respect to $B$.  (This is well-defined since $\bar Y$ and $B$ determine $\Sigma$.)

\medskip

\noindent{\bf Lemma~3.2.} {\sl Let $R$ be a complete d.v.r.\ and let $\bar Y$ be a projective normal connected $R$-curve, with general fibre $Y$.  Let $B$ be a proper closed 
subset of $Y$.  Let $(\bar Y_n, B, \Sigma_n)$ be the $n^{\rm th}$ blow-up of $(\bar Y, B, \Sigma)$, where $\Sigma$ is the intersection of the closure of $B$ in $\bar Y$ with the closed fibre of $\bar Y$.  Let $C$ be a proper closed subset of $Y$ that is disjoint from $B$.  Then for all sufficiently large $n$, the following conditions hold:

\cond{(i)} the closures of the points of $B$ are each regular $1$-dimensional subschemes of $\bar Y_n$ that meet the closed fibre of $\bar Y_n$ at distinct points;

\cond{(ii)} the closure of $C$ in $\bar Y_n$ is disjoint from $\Sigma_n$;

\cond{(iii)} the reduced closed fibre of $\bar Y_n$ is regular at the points of $\Sigma_n$.}

\medskip

\noindent{\it Proof.} For sufficiently large $n$, the closure of $B$ in $\bar Y_n$ is a regular $1$-dimensional scheme (of relative dimension $0$ over $R$), and so a disjoint union of irreducible regular components.  So (i) holds.  

Let $f$ be a rational function on $Y$ that vanishes on $C$ and has value $1$ on $B$.  Then for $n$ sufficiently large, $f$ is defined at the points of $\Sigma_n$ and 
the order of vanishing of $f$ is $0$ there.  So (ii) holds for such $n$.

By resolution of singularities for surfaces, after finitely many blow-ups the scheme $\bar Y$ will become regular.  So for sufficiently large $n$, the $n^{\rm th}$ blow-up $\bar Y_n$ of $\bar Y$ with respect to $B$ will be regular at each point $\sigma \in \Sigma_n$.  So the closure of $B$ is a Cartier divisor, given near $\sigma$ by some local defining function $f$.  Blowing up $\bar Y_n$ further with respect to $B$ will reduce the multiplicity of $f$ at the points of $\Sigma_m$ over $\sigma$  (for $m>n$) until it becomes $1$ at $\Sigma_N$, for some $N>n$.  At that point, the closure of $B$ will meet the closed fibre of $\bar Y_N$ over $\sigma$ only at a point of the last exceptional divisor, which is isomorphic to a projective line (since $\bar Y_m$ is regular over $\sigma$ for $m \ge n$, and in particular for $m=N-1$).  Hence after sufficiently many blow-ups, the reduced closed fibre will be regular where it meets the closure of $B$, giving (iii).   \qed

\medskip

\noindent{\bf Lemma~3.3.} {\sl Let $R$ be a complete d.v.r., let $G$ be a finite group, and suppose that 
$\psi:\bar Y \to \bar X$ is a $G$-Galois cover of projective normal connected $R$-curves, with general fibre $Y \to X$.  Suppose also that $B$ is the (reduced) inverse image of some proper closed subset $A \subset X$.  Let $\bar Y_n$ be the $n^{\rm th}$ blow-up of $\bar Y$ with respect to $B$.

a) Then for every $n \ge 0$, the $G$-action on $\bar Y$ lifts to a $G$-action on $\bar Y_n$, whose quotient $\bar X_n$ is equipped with a birational proper morphism to $\bar X$.  

b) If $n$ is sufficiently large, then distinct points of $A$ have the property that their closures in $\bar X_n$ meet the closed fibre of $\bar X_n$ at distinct points, and this closed fibre is locally irreducible at each of these points.}

\medskip

\noindent{\it Proof.}  (a) The set $B$ is $G$-invariant, and so the action of $G$ on $\bar Y$ extends inductively to each $\bar Y_n$.  The quotient $\bar X_n = \bar Y_n/G$ is proper over $R$ because $\bar Y_n$ is.  By the universal property of quotients, the composition $\bar Y_n \to \bar Y \to \bar X$ factors through $\bar X_n$; this gives a morphism $\bar X_n \to \bar X$ which is proper because $\bar X_n$ is proper over $R$.  This morphism is an isomorphism on the general fibre, because this property holds for the morphism $\bar X_n \to \bar X$.  So it is birational.  

(b) If $n$ is sufficiently large, then the closures of the points of $B$ in $\bar Y_n$ meet the closed fibre of $\bar Y_n$ at distinct points, by Lemma~3.2(i).  So the closures of the points of $A$ in $\bar X_n$ meet that closed fibre at distinct points (viz.\ the points in $\psi(\Sigma_n)$).  Also, by Lemma~3.2(iii), the closed fibre of $\bar Y_n$ is locally irreducible at the points of $\Sigma_n$.  But 
if there were distinct branches in the complete local ring of the closed fibre of $\bar X_n$ at some point $\xi$ in $\psi(\Sigma_n)$, then the closed fibre of $\bar Y_n$ would similarly have distinct branches locally at any point $\eta$ of $\Sigma_n$ over $\xi$ (since $\bar Y_n \to \bar X_n$ is a Galois branched cover).  So in fact the closed fibre of $\bar X_n$ is locally irreducible at each point of $\psi(\Sigma_n)$, i.e.\ at the points where the closure of $A$ meets the closed fibre of $\bar X_n$.   \qed

\medskip

Let $R$ be a complete d.v.r.\ with fraction field $\hat k$ and uniformizer $\pi$, let $\bar Y$ be a proper normal $R$-curve, and let $\eta, \eta'$ be closed points of the general fibre $Y$ of $\bar Y$.  
For a natural number $r$, we will say that the points $\eta, \eta'$ are {\it congruent modulo} $\pi^r$ if their closures in $\bar Y$ have the same pullbacks via $\Spec (R/\pi^rR) \to \Spec R$.  (Note that this depends on the model $\bar Y$ of $Y$, not just on the isomorphism class of the $\hat k$-curve $Y$.)

In particular, consider the projective $u$-line $\PP^1_R$, and a
closed point $\eta$ of the general fibre $\PP^1_{\hat k}$.  Suppose that the closure of $\eta$ in $\PP^1_R$ does not pass through the point at infinity.  Then $\eta$ is dominated by an $S$-point $(u=f)$ of $\PP^1_R$ for some finite extension $S$ of $R$ and some $f \in S$.  
If $g \in S$, then the closed point $\eta'$ that is dominated by $(u=g)$ will be congruent to $\eta$ modulo $\pi^r$ if $f,g$ are congruent modulo $\pi^r S$.  Since $\pi^r S$ is infinite, there are infinitely many possibilities for such $g$; so for every $r$ there are infinitely many points $\eta'$ with the same residue field that are congruent to $\eta$ modulo $\pi^r$.

\medskip

\noindent{\bf Lemma~3.4.} {\sl Let $k$ be a field, let $\hat k = k((t))$, and let $R = k[[t]]$.  Let $\Delta$ be a finite closed subset of the general fibre $\PP^1_{\hat k}$ of $\PP^1_R$.  

a) Let $P$ be a closed point 
of $\PP^1_R$ whose residue field is separable over $k$.  Then there is a closed point $\hat P$ of $\PP^1_{\hat k}$ whose residue field is separable over $\hat k$; whose closure in $\PP^1_R$ meets the closed fibre exactly at $P$; and such that no closed point of $\PP^1_{\hat k}$ that is congruent to $\hat P$ $({\sl mod}\, t^2)$ is contained in $\Delta$.  

b) Let $\bar Y \to \PP^1_R$ be a Galois branched cover of projective normal $R$-curves, with general fibre $Y \to \PP^1_{\hat k}$.  Let $\hat \ell$ be a finite extension of $\hat k$ that is contained in the function field of $Y$.  Assume that $\hat P'$ is a rational point of $\PP^1_{\hat \ell}$ that lies over a point $\hat P$ as in part (a).  Suppose moreover that $\hat P'$ splits completely under $Y \to \PP^1_{\hat \ell}$.  Let $h$ be a positive integer, and let $\Pi_h$ be
the set of rational points of $\PP^1_{\hat \ell}$ that are congruent to $\hat P'$ modulo $t^h$, are totally split under $Y \to \PP^1_{\hat \ell}$, have trivial decomposition group over $\PP^1_{\hat k}$, and have image in $\PP^1_{\hat k}$ that is not contained in $\Delta$.  Then the cardinality of $\Pi_h$ is equal to that of $\hat k$.}

\medskip

\noindent{\it Proof.}  (a) Let $k^s$ and $\hat k^s$ be the separable closures of $k$ and $\hat k$.  The point $P$ lies on the closed fibre $\PP^1_k \subset \PP^1_R$, which is the $u$-line over $k$.  Pick a closed point $\mu_0$ on $\PP^1_{k^s}$ over $P$; this is of the form $(u=c_0)$ for some $c_0 \in k^s$, and $P$ is given by the minimal polynomial of $c_0$ over $k$.

The closure $\bar \Delta$ of $\Delta$ in $\PP^1_R$ is a proper closed subset of $\PP^1_R$, and the reduction $\Delta_2$ of $\bar \Delta$ modulo $t^2$ is a proper closed subset of $\PP^1_R/(t^2)$.  Since $k^s$ is infinite, we may choose $c_1 \in k^s$ such that the locus of $(u=c_0 + c_1t)$ in $\PP^1_R/(t^2)$ is not contained in $\Delta_2$.  Let $c=c_0 + c_1t \in \hat k^s$ and let $\mu \in \PP^1_{\hat k^s}$ be the $\hat k^s$-point given by $(u=c)$.   Thus no $\hat k^s$-point of $\PP^1_{\hat k^s}$ that is congruent to $\mu$ modulo $t^2$ is contained in $\Delta$.  Also, $P$ is the closed point of $\PP^1_R$ where the closures of $\mu$ and its $\hat k$-conjugates each meet the closed fibre $\PP^1_k$.  So we may take $\hat P$ to be the closed point of $\PP^1_{\hat k}$ over which the  $\hat k^s$-point $\mu$ lies.

(b) Since $\hat k$ is infinite, and since $\hat \ell$ is a finite extension of $\hat k$, these two fields have the same cardinality, which is equal to that of ${\rm card}\,\PP^1(\hat \ell)$.  So ${\rm card}\,\Pi_h$ is bounded above by ${\rm card}\,\hat k$.
Also, $\Pi_h \supset \Pi_{h'}$ for $h<h'$.  So it suffices to show that ${\rm card}\,\Pi_h \ge {\rm card}\,\hat k$ for all sufficiently large $h$.

Let $S$ be the integral closure of $R$ in $\hat \ell$.  The closure of $\hat P'$ in $\PP^1_S$ meets the closed fibre at a unique point of $\PP^1_\ell$, where $\ell$ is the residue field of $S$.  After a change of variables, we may assume that this is not the point at infinity on $\PP^1_\ell$; and hence that $\hat P'$ is given by $(u=c)$ for some $c \in S$.  

Since the point $\hat P'$ splits completely, the strong form of Hensel's Lemma [La, II \S2, Prop.~2] implies that for all sufficiently large integers $h>1$ the cover $Y \to \PP^1_{\hat\ell}$ is also totally split over any rational point of $\PP^1_{\hat\ell}$ that is congruent modulo $t^h$ to $\hat P'$.  Since $\hat \ell$ is a finite separable extension of $\hat k$, there is a primitive element $d \in \hat \ell$ over $\hat k$.  For each non-zero element $a \in \hat k$, the element $ad \in \hat \ell$ is also a primitive element over $\hat k$, and the cover $Y \to \PP^1_{\hat\ell}$ is totally split over the $\hat\ell$-rational point $(u=c+adt^h)$ in $\PP^1_{\hat\ell}$.  Moreover this point has trivial decomposition group over its image in $\PP^1_{\hat k}$, because $ad$ is a primitive element for $\hat\ell$ over $\hat k$.  By part (a), this image cannot be contained in $\Delta$, since it is congruent to $\hat P$ modulo $t^h$ and $h>1$.  So this point lies in $\Pi_h$.  So ${\rm card}\,\Pi_h \ge {\rm card}\,\hat k$ for sufficiently large $h$,
as needed.  \qed

\medskip

The statement and proof of the next result are similar to those of [Ha87, Theorem~2.3].

\medskip

\noindent{\bf Proposition~3.5.} {\sl Let $k$ be a field of characteristic $p \ge 0$, and let $R = k[[t]]$.
Let $\Gamma$ be a finite group with a normal subgroup $N$, such that the quotient map $\Gamma \to G := \Gamma/N$ has a section $\sigma$.  Let $\tilde Y \to \tilde X$ be a $G$-Galois cover of normal proper $R$-curves, with general fibre $Y \to X$.  Let $n_1,\dots,n_r,m_1,\dots,m_s \ne 1$ be generators
for $N$ such that $p$ does not divide the order $o(n_i)$ of any $n_i$, and such that the order $o(m_j)$ of each $m_j$ is a power of $p$.  Suppose that the function field of $\tilde Y$ contains a primitive $o(n_i)^{\rm th}$ root of unity for each $i$.  

Let $\Phi_1,\dots,\Phi_r,\Psi_1,\dots,\Psi_s$ be rational functions on $\tilde X$ that restrict to the constant function $1$ at each generic point of the closed fibre, and whose zero and pole divisors are reduced.  Suppose that the supports of these functions on $\tilde X$ are pairwise disjoint and contained in the smooth locus over $R$; that $\tilde Y \to \tilde X$ splits completely over each closed point in these supports; and that their supports are disjoint from 
some proper closed subset $D \subset Y$ containing the ramification locus of $Y \to X$.  
Then there is an irreducible normal $N$-Galois cover $\tilde Z \to \tilde Y$ such that the following conditions hold:

\cond{(i)} $\tilde Z \to \tilde Y$ is branched precisely over the supports of the divisors of the $\Phi_i$'s and the zero divisors of the $\Psi_j$'s, with ramification indices $o(n_i)$ and $o(m_j)$ respectively;

\cond{(ii)} the composition $\tilde Z \to \tilde X$ is $\Gamma$-Galois;

\cond{(iii)} the field $k$ has the same algebraic closures in the function fields of $\tilde Y$ and $\tilde Z$;

\cond{(iv)} the closed fibre of $\tilde Z \to \tilde Y$ is an $N$-Galois mock cover; {\new and for each $1 \le i \le r$ [resp.\ $1 \le j \le s$] there is a closed point $Q$ {\new on} the identity sheet {\new of the closed fibre of $\tilde Z$} lying over a closed point $P$ in the closure of the divisor of $\Phi_i$ [resp.\ the zero divisor of $\Psi_j$] such that $n_i$ [resp.\ $m_j$] generates the inertia group {\new at} $Q$ and {\new at} the generic point of each of the irreducible components of the divisor {\new on $\tilde Z$ defined by} $\Phi_i$ [resp.\ the zero divisor {\new on $\tilde Z$ defined by} $\Psi_j$] that passes through $Q$.}

\cond{(v)} the points $\delta$ of $D$ are totally split under $\tilde Z \to \tilde Y$;

\cond{(vi)} the decomposition groups of $\tilde Z \to \tilde X$ at the points of $\tilde Z$ over any $\delta \in D \subset \tilde Y$ are the conjugates of $\sigma(G_\delta)$, where $G_\delta$ is the decomposition group of $\tilde Y \to \tilde X$ at $\delta$.

} 

\medskip

\noindent{\it Proof.}  For each $i=1,\dots,r$, and each closed point $P \in \tilde X$ in the support of the divisor of $\Phi_i$, pick a closed point $Q \in \tilde Y$ over $P$.  Thus each point in the fibre of $\tilde Y \to \tilde X$ over $P$ is of the form $g(Q)$ for some unique $g \in G$, because $\tilde Y \to \tilde X$ is totally split over $P$.  

For each $g \in G$, let ${\cal L}_{g(Q)}$ be the fraction field of the complete local ring $\hat {\cal O}_{\tilde Y,g(Q)}$.  Also, let
${\cal A}_{g(Q)}$ be the normalization of $\hat {\cal O}_{\tilde Y,g(Q)}$ in the field extension of ${\cal L}_{g(Q)}$ given by $z^{o(n_i)}=\Phi_i$.   
Since $\zeta_{o(n_i)}$ is contained in the function field of $\tilde Y$ and hence in ${\cal L}_{g(Q)}$, this is a cyclic extension of degree $o(n_i)${\new, totally ramified over $g(Q)$ and the divisor of $\Phi_i$}.  Moreover, its Galois group may be identified with $gN_ig^{-1}$, where $N_i$ is the subgroup of $N$ generated by $n_i$.  Since $\Phi_i$ restricts to the constant function $1$ at the generic point $g(Q)^\circ$ of the closed fibre of $\Spec \hat {\cal O}_{\tilde Y,g(Q)}$, the pullback of the 
cyclic cover $\Spec {\cal A}_{g(Q)} \to \Spec \hat {\cal O}_{\tilde Y,g(Q)}$ to the complete local ring $\hat {\cal O}_{g(Q)^\circ}$ at this generic point is trivial.  
This cover is equipped with an induced indexing, by the elements of $gN_ig^{-1}$, of the components of the restriction to $\Spec \hat {\cal O}_{g(Q)^\circ}$.  

Taking a disjoint union of copies of this cover, indexed by the cosets of $gN_ig^{-1}$ in $N$, we obtain an induced (disconnected) $N$-Galois cover of $\Spec {\cal A}_{g(Q)}$, together with an indexing by $N$ of the components over $\hat {\cal O}_{g(Q)^\circ}$.  As {\new $g$ varies over $G$}, 
these indexings are compatible with the action of $G$ on $\tilde Y$; and so we obtain a (disconnected) $\Gamma$-Galois cover of $\Spec \hat {\cal O}_{\tilde X,P}$,
viz.\ $\Ind_{N_i}^\Gamma \Spec {\cal A}_Q$,
which is ramified precisely over the support of the divisor of $\Phi_i$ on $\Spec \hat {\cal O}_{\tilde X,P}$
(using that the zero and pole loci of $\Phi_i$ are reduced). {\new  Taking {now $g={\rm id}_G$}, we see that $n_i$ generates the {\new inertia group of the components of the inverse image of this divisor} that meet the identity sheet of the closed fibre.}

Similarly, for each $j=1,\dots,s$, pick a closed point $Q \in \tilde Y$ over each closed point $P \in \tilde X$ in the support of the zero divisor of $\Psi_j$, and again consider the field ${\cal L}_{g(Q)}$ for each $g \in G$.  We have that $o(m_j) = p^{\beta_j}$, for some positive integer $\beta_j$.  In the $\beta_j^{\rm th}$ truncated Witt vector ring $W_{\beta_j}(\hat {\cal O}_{\tilde Y,g(Q)}[y_1,\dots,y_{\beta_j}])$, let $\underline\Psi_j$ and $\underline y$ denote the elements with Witt coordinates $(\Psi_j,0,\dots,0)$ and $(y_1,\dots,y_{\beta_j})$ respectively, and let ${\rm Fr}$ denote Frobenius.  Consider the field extension of ${\cal L}_{g(Q)}$ given by the Witt coordinates of ${\rm Fr}(\underline y)-\underline\Psi_j^{p-1}\underline y = t$, and let ${\cal A}_{g(Q)}$ be the normalization of $\hat {\cal O}_{\tilde Y,g(Q)}$ in this extension.  
This is a cyclic extension whose Galois group may be identified with $gN_j'g^{-1}$, where $N_j'$ is the subgroup of $N$ generated by $m_j$.  Proceeding as before, we obtain an induced (disconnected) $N$-Galois cover of $\Spec {\cal A}_{g(Q)}$, whose restriction to $\hat {\cal O}_{g(Q)^\circ}$ is trivial, with components indexed by the elements of $N$.  Letting $g$ vary and taking the union, we obtain a $\Gamma$-Galois cover of $\Spec \hat {\cal O}_{\tilde X,P}$,
viz.\ $\Ind_{N_j'}^\Gamma \Spec {\cal A}_Q$,
which is ramified precisely over the zero locus of $\Psi_j$ on $\Spec \hat {\cal O}_{\tilde X,P}$ (using that this locus is reduced){\new, and with $m_j$ generating the inertia group on the components of this locus that meet the identity sheet of the closed fibre}.

Let $\tilde Y_0$ be the closed fibre of $\tilde Y$; let $\Delta$ be the set of points where the support of some $\Phi_i$ or $\Psi_j$ meets $\tilde Y_0$; let $U$ be the complement of $\Delta$ in $\tilde Y_0$; and let ${\cal U}$ be the completion of $\tilde Y$ along $U$.  (That is, ${\cal U} = \Spec \O_{\cal U}$, where $\O_{\cal U}$ is the $t$-adic completion of the ring of functions on an affine open subset $\tilde U \subset \tilde Y$ such that $\tilde U \cap \tilde Y_0 = U$.)  Consider the trivial $N$-Galois cover of ${\cal U}$, consisting of $|N|$ disjoint copies of ${\cal U}$ that are indexed by the elements of $N$, and on which $N$ acts by the regular representation.  By formal patching (e.g.\ [HS99, \S1, Cor.\ to Thm.~1], [Pr, Thm.~3.4]), there is an $N$-Galois cover $\tilde Z \to \tilde Y$ whose pullback to the spectrum of each $\hat {\cal O}_{\tilde Y,g(Q)} = g(\hat{\cal O}_{\tilde Y,Q})$ is isomorphic to the cover described above; whose pullback to ${\cal U}$ is trivial; and such that over each $g(Q)$, the two isomorphisms with the trivial cover (induced from pulling back the isomorphism over ${\cal U}$ and the one over $\hat {\cal O}_{\tilde Y,g(Q)}$) give the same indexing of components.  Since these indexings were compatible (for various $g \in G$), the action of $G = {\rm Gal}(\tilde Y/\tilde X)$ lifts to an action on $\tilde Z$, compatibly with the conjugation action of $G$ on $N$ in $\Gamma$; and so the composition $\tilde Z \to \tilde X$ is Galois with group $\Gamma$.   Since the local $N$-covers given above are normal, so is $\tilde Z$ (as normality is a local property).

The closed fibre of $\tilde Z \to \tilde Y$ is an $N$-Galois mock cover, 
consisting of a union of copies of $\tilde Y_0$, indexed by the elements of $N$.  For each point $P \in \tilde X$ in the closure of the support of the divisor of some $\Phi_i$ [resp.\ the zero divisor of some $\Psi_j$], consider the chosen point $Q \in \tilde Y$ over $P$.  The closed fibre $\tilde Z_0$ of $\tilde Z$ is a union of (intersecting) copies of $\tilde Y_0$ indexed by $N$; let $\tilde Q \in \tilde Z_0$ be the closed point on the identity copy of $\tilde Y_0$ corresponding to $Q$.  
The inertia group of $\tilde Z\to \tilde Y$ at the point $\tilde Q$ is generated by $n_i$ [resp.\ $m_j$].  Since the $n_i$'s and $m_j$'s generate $N$, it follows that the closed fibre of $\tilde Z$ is connected.  Hence $\tilde Z$ is connected and thus irreducible (since $\tilde Z$ is normal).  
Moreover, if $\hat \ell$ is the algebraic closure of $\hat k$ in the function field of $\tilde Y$, then 
$\hat\ell$ is algebraically closed in the function field of $\tilde Z$, since the closed fibre $\tilde Z_0 \to \tilde Y_0$ is a mock cover.  

It remains to verify (v) and (vi).  The cover $\tilde Z \to \tilde Y$ is totally split over each point $\delta$ of $D$, since those points are in the image of ${\cal U} \to \tilde Y$ and since the pullback of $\tilde Z \to \tilde Y$ to $\cal U$ is trivial.  Also, since that pullback is trivial, it consists of a disjoint union of copies of $\cal U$, indexed by the cosets of $N$ in $\Gamma$, and acted upon by $\Gamma$.  The identity copy and $\cal U$ itself are isomorphic, together with their $G$-actions (if we identify $G$ with $\sigma(G) \subset \Gamma$).  The other copies are isomorphic to $\cal U$ as schemes, but their group actions are conjugated by the corresponding elements of $N$.  So if we denote by $\delta_n$ the point of $\tilde Z$ over $\delta$ on the copy of $\cal U$ indexed by $n \in N$, then the decomposition group of $\tilde Z \to \tilde X$ at the point $\delta_n$ is the conjugate of $\sigma(G_\delta)$ by $n$.  Since $\Gamma$ is generated by $N$ and $\sigma(G)$, the decomposition groups of $\tilde Z \to \tilde X$ at the points of $\tilde Z$ over $\delta \in \tilde Y$ are precisely the $\Gamma$-conjugates of $G_\delta$.  So (v) and (vi) also hold.
\qed

\bigskip

\noindent {\bf Section~4. Split embedding problems over curves}

\medskip

In this section we extend Pop's result on solving split embedding problems for curves over large fields.  That extension (Theorem~4.1 below) will be used in next section in the proof of Theorem~1.1 (which appears there as Theorem~5.1).

Pop's result [Po96, Main Theorem~A] showed that such embedding problems have proper regular solutions; and a later result of Haran and Jarden [HJ, Theorem~6.4] showed that solutions can be chosen so as to split completely over a specified unramified point of the given cover.  (Actually, the results proven in [Po96] and [HJ] and in unpublished manuscripts of those authors were somewhat less general than this; see [Ha03, Remark~5.1.11].)  

Here we show that more is true:  that the solutions to the embedding problem can be chosen to be totally split over any finite set of (possibly ramified) points; that there are ``many'' such solutions; and for a large field of the form $k((t))$, that the solution to the embedding problem can be chosen to be totally ramified over a specified closed point on a model over $k[[t]]$.  These assertions are contained in Theorem~4.1 in the case of $k((t))$.  Afterwards, in Theorem~4.3, we obtain a more general result for curves over large fields that need not be of the form $k((t))$ (though at the expense of no longer being able to speak of an integral model).  Then in Corollary~4.4, we deduce a result about absolute Galois groups of function fields over ``very large'' fields being {\new semi-free}.

\medskip

\noindent {\bf Theorem~4.1.} {\sl Let $k$ be an arbitrary field and let $X$ be a smooth connected projective curve over $\hat k = k((t))$, with function field $K$.

a) Then every finite split embedding problem ${\cal E} = (\alpha:G_K \to G, f:\Gamma \to G)$ for $K$ has a proper regular solution $\beta:G_K \to \Gamma$.  

b) For such an $\E$, with section $\sigma$ of $f$, 
let $\pi:Y \to X$ be the $G$-Galois branched cover of normal curves corresponding to $\alpha$, and let $D \subset Y$ be a finite set of closed points.  Then we may choose the proper regular solution in (a) so that the corresponding normal cover $Z \to Y$ is totally split over the points of $D${\new ;}
the decomposition groups of $Z \to X$ at the points of $Z$ over $\delta \in D$ are the conjugates of $\sigma(G_\delta)$, where $G_\delta$ is the decomposition group of $Y \to X$ at $\delta${\new; and $Z \to Y$ has no non-trivial \'etale subcover $W \to Y$}.

c) If $\E$ is non-trivial then there are  {\new exactly} ${\rm card}(\hat k)$ {\new non-isomorphic}
proper regular solutions $Z \to Y \to X$ to $\cal E$ satisfying (b){\new, and there is a set of ${\rm card}(\hat k)$ such solutions that are {\new independent}}.

d) Suppose that $\bar X$ is a smooth projective model for $X$ over $R{\new = k[[t]]}$; and let $\bar Y$ be the normalization of $\bar X$ in $Y$.  Suppose also that $Q$ is a closed point of $\bar Y$ such that the residue field of its image $P$ in $\bar X$ is separable over $k$.  If $\E$ is non-trivial then there are ${\rm card}(\hat k)$ 
{\new independent} proper regular solutions $Z \to Y \to X$ in (c) such that the normalization $\bar Z$ of $\bar Y$ in $Z$ is totally ramified over the closed point $Q \in \bar Y$.  The pullbacks $Z^* \to Y^*$ via $Y^* := {\rm Spec}\,{\hat {\cal O}}_{{\bar Y},Q} \to {\bar Y}$ of the corresponding $\ker(f)$-Galois covers $\bar Z \to \bar Y$ are each irreducible and 
{\new form a set of linearly disjoint covers of} cardinality equal to ${\rm card}(\hat k)$.  Moreover, if $Q$ is the only closed point of $\bar Y$ over $P \in\bar X$, then the pullbacks of the covers $\bar Z \to \bar X$ under $X^* := {\rm Spec}\,{\hat {\cal O}}_{{\bar X},P} \to {\bar X}$ form a set of irreducible 
$\Gamma$-Galois covers $Z^* \to X^*$ having the same cardinality.

{\new e) Suppose in (d) that $\bar X = X_0 \times_k R$, for some smooth connected projective $k$-curve $X_0$.  Then the asserted proper solutions in (d) may be chosen so that $Z \to Y$ is \'etale over every point of $Y$ whose image in $X$ is of the form $P' \times_k \hat k$ with $P'$ a point of $X_0$.}}

\medskip

In using Theorem~4.1 to prove Theorem~5.1, we will start with a $G$-Galois cover of $\Spec k[[x,t]]$; extend it to a cover of $x$-line over $k[[t]]$; use Theorem~4.1 to dominate that by a $\Gamma$-Galois cover; and then restrict that back to $\Spec k[[x,t]]$.  The condition on total ramification over a closed point (in part (d) of Theorem~4.1 above) is then used to conclude that the restricted cover remains connected.  In the proof of Theorem~4.1, the condition on total ramification is obtained by blowing up and down in such a way that the branch components are all forced to pass through the designated point.  The assertion in part (b) of Theorem 4.1 strengthens the splitting condition of Haran-Jarden; and the ``many {\new independent} solutions'' assertion 
{\new in parts (c) and (d)} 
allows us to conclude in Section~5 that the absolute Galois group of $k((x,t))$ is {\new semi-free}.

Theorem~4.1 is proved with the aid of Proposition~4.2, which reduces the theorem to a more manageable special case (and whose proof uses results from Section~3).  

\medskip

\noindent {\bf Proposition~4.2.} {\sl Let $p = {\rm char}\ k$.  In order to prove Theorem~4.1, it suffices to prove it in the special case that the $G$-Galois cover $\pi:Y \to X$ corresponding to $\alpha$ has the following properties:

(i) the function field of
$Y$ contains a primitive $d^{\rm th}$ root of unity, where $d$ is the maximal prime-to-$p$ factor of $|\ker f|$;

(ii) there is a finite generically unramified $\hat k$-morphism $\phi:X \to 
\PP^1_{\hat k}$ such that the composition $\phi \pi:Y \to \PP^1_{\hat k}$ is Galois{\new;}  
in part (d) also such that $\phi$ extends to a finite $R$-morphism $\bar X \to \PP^1_R$ (where in parts (a) - (c), $\bar X$ denotes the normalization of $\PP^1_R$ in $X$, viewing $\PP^1_{\hat k}$ as the generic fibre of $\PP^1_R$); {\new and in part (e) also such that $\phi:\bar X \to \PP^1_R$ is the base change from $k$ to $R$ of a finite $k$-morphism $\phi_0:X_0 \to \PP^1_k$;}

(iii) there is a closed point $\hat P \in X$ such that the points in the fibre of $Y \to X$ over $\hat P$ are rational over the algebraic closure of $\hat k$ in the function field of $Y$; no closed point of $X$ that is congruent to $\hat P$ modulo $t^2$ (with respect to $\bar X$) is contained in the branch locus of $Y \to X$ or in the image of $D$; $\hat P$ specializes at the closed fibre to a closed point $P \in \bar X$ whose residue field is separable over $k$; and in part (d), $P$ is the image of the given point $Q \in \bar Y$.}

\medskip

\noindent{\it Proof.}  Suppose that Theorem~4.1 is proven in the above special case.  For part (d) of Theorem~4.1, this special case includes the additional assumptions that $\phi$ extends to the given model $\bar X$, and that the specialization $P$ of $\hat P$ is the image of the given point $Q$.  

Now assume that we are in the general case of Theorem~4.1, which we wish to prove.  Let $N = \ker f$.  In part (b), after enlarging $D$, we may assume that $D$ contains the ramification locus of $Y \to X$.  (In part (a), we let $D$ denote the ramification locus of $Y \to X$.)

According to Proposition~3.1(a), there is a finite generically unramified $R$-morphism $\bar X \to \PP^1_R$ for some proper normal model $\bar X$ for $X$ over $R$.  By Proposition~3.1(b), 
in part (d) we may choose $\bar X \to \PP^1_R$ with respect to the given $R$-model $\bar X$ of $X$. {\new In part (e), we may choose a finite generically unramified $k$-morphism $\phi_0:X_0 \to \PP^1_k$; and we may then choose $\bar X \to \PP^1_R$ to be the base change of $\phi_0$ from $k$ to $R$.}

Let $P_0 \in \PP^1_R$ be a closed point whose residue field is separable over $k$; in part (d) we may choose this point to be the image of $Q \in \bar Y$ in $\PP^1_R$ (since the residue field at $Q$ is separable over $k$).  Also let $\Delta$ be the 
union of the branch locus of $X \to \PP^1_{\hat k}$ with the
image of $D$ in $\PP^1_{\hat k}$.  Applying Lemma~3.4(a) to $P_0$ and $\Delta$, we obtain a closed point $\hat P_0$ of $\PP^1_{\hat k}$ whose residue field is separable over $\hat k$; whose closure in $\PP^1_R$ meets the closed fibre exactly at $P_0$; and such that no closed point of $\PP^1_{\hat k}$ that is congruent to $\hat P_0$ modulo $t^2$ (with respect to $\PP^1_R$) is contained in $\Delta$.  
In particular, $\hat P_0$ is not in the branch locus of $Y \to \PP^1_{\hat k}$; hence the residue fields at the points of $X$ and $Y$ over $\hat P_0$ are separable over $\hat k$.

Let $\tilde Y \to \PP^1_{\hat k}$ be the Galois closure of $Y \to \PP^1_{\hat k}$.  So its function field $\tilde L$ is Galois over $\hat k(u)$, the function field of $\PP^1_{\hat k}$.  Let $\hat\ell$ be the algebraic closure of $\hat k$ in $\tilde L$.  Let $\hat\ell'$ be a finite field extension of $\hat\ell$ that is Galois over $\hat k$ and that contains a primitive $d^{\rm th}$ root of unity, and such that every closed point of $Y$ lying over $\hat P_0$ has residue field contained in $\hat\ell'$. 

Let $Y' = \tilde Y \times_{\hat\ell} \hat\ell'$, with function field equal to the compositum $\tilde L\hat\ell'$ in a separable closure of $K$.  Since $\tilde L/\hat k(u)$ and $\hat\ell'/\hat k$ are Galois, so is $\tilde L\hat\ell'/\hat k(u)$; i.e.\ $Y'
\to \PP^1_{\hat k}$ is Galois.  Moreover this morphism extends to a finite morphism $\bar Y' \to \PP^1_R$, where $\bar Y'$ is the normalization of $\bar X$ in $Y'$.  Note that $\bar Y' \to \bar X$ factors through $\bar Y$, the normalization of $\bar X$ in $Y$.

The intermediate cover 
$Y' \to X$ is necessarily also Galois, with Galois group 
$G' := {\rm Gal}(\tilde L\hat\ell'/K)$; this is an extension of $G = {\rm Gal}(L/K)$ through which $\alpha:G_K \to G$ factors, via an epimorphism $\alpha':G_K \to G'$.  Let $\E_{\alpha'} = (\alpha':G_K \to G', f':\Gamma' \to G')$ be the induced embedding problem, where $\Gamma' =  \Gamma \times_G G'$ (see Section~2).  This is also a finite embedding problem, and it is split by Lemma~2.6, in the context of Example~2.5(b) there.  Also, $\ker f' = N \times 1 \approx N = \ker f$.  So the $G'$-Galois cover $Y' \to X$ satisfies (i) and (ii), with respect to the embedding problem $\E_{\alpha'}$, which has the same kernel as $\E$.

Choose a closed point $P \in \bar X$ lying over $P_0 \in \PP^1_R$; in part (d) we may take $P$ to be the image of $Q \in \bar Y$.  Let $\hat P$ be a closed point of $X$ that specializes to $P$ and that lies over $\hat P_0$; this exists by the going-down theorem applied to $\bar X \to \PP^1_R$.  So $\hat P$ satisfies (iii).

Thus the finite split embedding problem $\E_{\alpha'} = (\alpha':G_K \to G', f':\Gamma' \to G')$, with $\alpha'$ corresponding to $Y' \to X$, satisfies the three conditions of Proposition~4.2.  So by assumption, the conclusion of Theorem~4.1 holds for this embedding problem.  

By (a) of Theorem~4.1, there is a proper regular solution to $\E_{\alpha'}$, corresponding to an $N$-Galois cover $Z' \to Y'$.  By Lemma~2.4 and Proposition~2.9, every proper regular solution to $\E_{\alpha'}$ yields a proper regular solution to $\E$, via the map ${\rm PS}(\E_{\alpha'}) \to {\rm PS}(\E)$.  So we obtain a proper regular solution to $\E$, corresponding to a cover $Z \to Y$.  This proves (a) for $\E$. 

For (b), let $D' \subset Y'$ be the inverse image of $D \subset Y$, and let
$Z' \to Y'$ correspond to a proper regular solution to $\E_{\alpha'}$ over which $D'$ splits completely, which exists by the special case.  Then the decomposition groups of $Z' \to X$ at the points of $Z'$ over $\delta' \in D'$ are the conjugates of $\sigma'(G'_{\delta'})$, where $G'_{\delta'}$ is the decomposition group of $Y' \to X$ at $\delta'$ and where $\sigma'$ is the splitting of $f'$ given by Lemma~2.6.  So 
for the cover $Z \to Y$ corresponding to the proper regular solution obtained for $\E$, the decomposition groups of $Z \to X$ at the points of $Z$ over $\delta \in D$ are the conjugates of $\sigma(G_\delta)$.  So the decomposition groups of $Z \to X$ at these points are mapped isomorphically onto that of $Y \to X$ at $\delta$; hence $Z \to Y$ is totally split over $\delta$.  
{\new Moreover, if $W \to Y$ is an \'etale subcover of $Z \to Y$, then
$W' := W \times_Y Y' \to Y'$ is an \'etale subcover of $Z' \to Y'$.  Since the map $\ker f \to \ker f'$ is an isomorphism, $\Gal(Z'/W') = \Gal(Z/W)$, and $W \to Y$ is non-trivial if and only if $W' \to Y'$ is.  Since $W' \to Y'$ is trivial by part (b) in the special case, $W \to Y$ is also {\new trivial.}}
So (b) holds for the given embedding problem $\E$.

For (c), by the special case applied to $\E_{\alpha'}$, we have that $\E_{\alpha'}$ has ${\rm card}(\hat k)$ {\new independent} solutions to that embedding problem satisfying (a) and (b) there.  But by Lemma~2.4, the map ${\rm PS}(\E_{\alpha'}) \to {\rm PS}(\E)$ is injective{\new, and moreover it maps a set of independent proper solutions to a set of independent proper solutions.}  So {\new the given solutions to $\E_{\alpha'}$ induce a set of}  
${\rm card}(\hat k)$ {\new independent} solutions to $\E$ above, corresponding to connected covers $Z \to Y\to X$ of $\hat k$-curves {\new that are linearly disjoint over $Y$}, and to irreducible normal $R$-models $\bar Z \to \bar Y \to \bar X$.
But the cardinality of the function field of $X$ is equal to ${\rm card}(\hat k)$; so in fact the cardinality of 
{\new the set of proper regular solutions to $\E$}
is at most ${\rm card}(\hat k)$ and hence exactly ${\rm card}(\hat k)$.   
This proves (c).

For (d), let $Q' \in \bar Y'$ be a closed point over $Q \in \bar Y$.  Then $Q,Q'$ have the same image $P$ in $\bar X$, which has separable residue field.  So by part (d) of the special case of Theorem~4.1, the proper regular solution $Z' \to X'$ to $\E_{\alpha'}$ can be chosen, in ${\rm card}(\hat k)$ 
{\new independent} 
ways, so that the normalization $\bar Z'$ of $\bar X$ in $Z'$ is totally ramified over $Q'$.   
For each of these solutions, the above cover $\bar Z \to \bar Y$ is totally ramified over $Q \in \bar Y$, since $\bar Z' \to \bar Y'$ is the normalized pullback of $\bar Z \to \bar Y$ by Lemma~2.9(b).
So the first part of (d) holds, again using the injectivity of  ${\rm PS}(\E_{\alpha'}) \to {\rm PS}(\E)$ {\new and of ${\rm PS}(\E^n_{\alpha'}) \to {\rm PS}(\E^n)$}.   

Next, we consider the cardinality of the pullbacks $Z^* \to Y^*$ of $\bar Z \to \bar Y$.  First, by this part of (d) in the above special case, there {\new is a set of} ${\rm card}(\hat k)$ {\new linearly disjoint} pullbacks $Z'^* \to Y'^*$ of the above solutions $\bar Z' \to \bar Y'$ for $\E_{\alpha'}$, where $Y'^* = \Spec 
{\hat {\cal O}}_{{\bar Y'},Q'}$ and where $Q' \in \bar Y'$ is a point over $Q \in \bar Y$.  Let ${\cal F}$ be the embedding problem $(G_{L^*}\to 1, N \to 1)$, where $L^*$ is the function field of ${\hat {\cal O}}_{{\bar Y},Q}$; let $J$ be the Galois group of the function field of $Y'^*$ over that of $Y^*$; 
let $\beta:G_{L^*} \to J$ be the corresponding surjection;
and let ${\cal F}_\beta$ be the induced embedding problem with respect to $\beta:G_{L^*} \to J$.  (So ${\cal F}, \beta, J$ here play the roles of the objects ${\cal E}, \alpha', G'$ in the definition of induced embedding problem in Section 2.)
The above pullbacks $Z'^* \to Y'^*$ form a subset of ${\rm PS}({\cal F}_{\beta})$; and their images $Z^* \to Y^*$ under the injection ${\rm PS}({\cal F}_{\beta}) \to {\rm PS}({\cal F})$ are precisely the pullbacks to $Y^*$ of the above solutions $\bar Z \to \bar Y$ to $\E$.  So the cardinality of these pullbacks $Z^* \to Y^*$ is 
at least ${\rm card}(\hat k)$.  Since the cardinality of the function field of $Y^*$ is equal to ${\rm card}(\hat k)$, this inequality is actually an equality.   {\new Moreover, these ${\rm card}(\hat k)$ solutions $Z^* \to Y^*$ are {\new independent}, because the solutions $\bar Z'^* \to \bar Y'^*$ are, using the injectivity of the map  
${\rm PS}({\cal F}^n_{\beta}) \to {\rm PS}({\cal F}^n)$ for each $n$ (where these latter embedding problems are defined via $n$th fibre powers as in the case of $\E^n_{\alpha'}$ and $\E^n$).}
So this part of (d) is verified.

{\new To conclude the proof of (d),} if two $N$-Galois covers $Z^* \to Y^*$ are non-isomorphic, then their compositions with $Y^* \to X^*$ are also non-isomorphic as (possibly disconnected) $\Gamma$-Galois covers.  But these compositions are irreducible $\Gamma$-Galois covers if $Q$ is the unique point of $\bar Y$ over $P \in \bar X$, so the last part of (d) follows.  {\new 

Finally, part (e) holds for the given embedding problem $\E$ because the branch locus of $Z \to Y$ is contained in that of $Z' \to Y'$.}
\qed

\medskip

Using the above reduction result, we now prove Theorem~4.1.  That is, given a finite $G$-Galois cover $Y \to X$, where $G = \Gamma/N$, we want to construct a $\Gamma$-Galois cover $Z \to X$ that dominates the given $G$-Galois cover and has certain additional properties.  The proof will use Proposition~3.5; and for this, we need to construct rational functions $\Phi_i, \Psi_j$ corresponding to generators of $G$, as in the proposition.  Using Proposition~4.2, it will be sufficient to construct these rational functions on $\PP^1_{\hat k}$ rather than on $X$ (which maps to $\PP^1_{\hat k}$).  These functions will be obtained as norms (from $\hat \ell$ to $\hat k$, where $\hat \ell$ is the algebraic closure of $\hat k$ in the function field of $Y$) of rational functions $\varphi_i, \psi_j$ on $\PP^1_{\hat \ell}$.  The divisors of these functions will have the form $\hat P_i' - \hat P_i'^*$ and $\hat P_j'' - \hat P_j''^*$ respectively, for some $\hat\ell$-points $\hat P_i',\hat P_i'^*,\hat P_j'',\hat P_j''^*$ on the projective line over which $Y \to \PP^1_{\hat \ell}$ is totally split.  These points are constructed with the aid of Lemma 3.4(b), and then yield the functions $\varphi_i, \psi_j$ and $\Phi_i, \Psi_j$.

\medskip

\noindent{\it Proof of Theorem~4.1.}  {\new We may assume that the embedding problem is non-trivial, since otherwise the result is immediate.} 
We proceed in several steps, following the strategy outlined above.

\smallskip

\noindent {\it Step 1.}  Construction of the points $\hat P_i',\hat P_j''$:

\smallskip
 
By Proposition~4.2, we may assume that the hypothesis (i)-(iii) there hold; and we preserve the notation used there.  So we have a morphism $\phi:X \to 
\PP^1_{\hat k}$ and a closed point $\hat P \in X$ specializing to a point $P$ on the closed fibre, having the properties stated in Proposition~4.2.  In particular, in part (d) of the theorem, $P$ is the image of the given closed point $Q \in \bar Y$.  Let $E$ be the Galois group of the cover $Y \to \PP^1_{\hat k}$ given by (ii) of Proposition~4.2.  Note that $Y \to \PP^1_{\hat k}$ is the generic fibre of $\bar Y \to\PP^1_R$, and that the latter cover is also Galois with group $E$.

Let $N$ be the kernel of the surjection $\Gamma \to G$.   We may choose non-trivial generators $n_1,\dots,n_r, \allowbreak m_1,\dots,m_s$ for $N$ such that $p = {\rm char}\,k$ does not divide the order $o(n_i)$ of any $n_i$, and such that the order $o(m_j)$ of each $m_j$ is a power of $p$.  Here $r,s \ge 0$.  {\new Since we are assuming that $\E$ is non-trivial, the kernel $N$ is} 
non-trivial and so $r,s$ are not both $0$.  
By hypothesis (i) of Proposition~4.2, the function field of $Y$ contains a primitive $o(n_i)^{\rm th}$ root of unity for each $i$.

Let $\hat \ell$ be the algebraic closure of $\hat k$ in the function field of $Y$ and let $S$ be the integral closure of $R$ in $\hat \ell$.  So $\hat \ell$ is separable over $\hat k$ since $Y \to X$ is generically separable.  By (iii) of Proposition~4.2, the points in the fibre of $Y \to X$ over $\hat P \in X$ are $\hat \ell$-rational.  Pick such a point $\hat Q \in Y$ in this fibre; in part (d) we may choose it so that it specializes to the given point $Q$ on the closed fibre of $\bar Y$.  (In parts (a)-(c)), we choose $\hat Q$ arbitrarily in this fibre, and we let $Q$ denote the specialization of $\hat Q$ to the closed fibre of $\bar Y$.  Its residue field need not be separable over $k$.)

Let $P_0$ be the image of $P$ in $\PP^1_R$ and let 
$P'_0$ be the image of $Q$ in $\PP^1_S$.  So $P'_0$ lies over $P_0$.
Similarly, let $\hat P_0$ be the image of $\hat P$ in $\PP^1_{\hat k}$ and let 
$\hat P'_0$ be the image of $\hat Q$ in $\PP^1_{\hat \ell}$.  Thus 
$\hat P'_0$ is a rational point of $\PP^1_{\hat \ell}$ that lies over $\hat P_0$.  Let $\Delta \subset \PP^1_{\hat k}$ be a finite closed subset containing the branch locus of the composition $Y \to X \to \PP^1_{\hat k}$, together with the image of $D$ under this composition (in parts (b)-(d)).  Choose an integer $h>1$.  By Lemma~3.4(b), there are 
{\new ${\rm card}(\hat k)$} rational points of $\PP^1_{\hat \ell}$ that are congruent to $\hat P'_0$ modulo $t^h$, are totally split under $Y \to \PP^1_{\hat \ell}$, have trivial decomposition group over $\PP^1_{\hat k}$, and have image in $\PP^1_{\hat k}$ that is not contained in $\Delta$ (in fact, there are as many as the cardinality of $\hat k$).  Pick such points $\hat P'_i$ and $\hat P''_j$
for $1 \le i \le r$ and $1 \le j \le s$, no two of which are conjugate over $\hat k$.  
{\new (There are ${\rm card}(\hat k)$ ways of choosing this set of points; and in particular, in part (e), we may choose them not to include the point $P_0' \times_\ell \hat \ell$.)}
The closures of these points, and their $\hat k$-conjugates, meet the closed fibre of $\PP^1_S$ at $P'_0$ and its $k$-conjugates.
 
\smallskip

\noindent{\it Step 2.}  Blowing up to separate the closures of the points $\hat P_i',\hat P_j''$:

\smallskip

Let $A_0 \subset \PP^1_{\hat\ell}$ be the set consisting of all the $\hat\ell$-points $\hat P'_i$ and $\hat P''_j$, and let $A \subset \PP^1_{\hat\ell}$ consist of the elements of $A_0$ and their $\hat k$-conjugates.  Let $B_0,B \subset Y$ be the inverse images of $A_0,A$ respectively; these are sets of $\hat\ell$-points of $Y$, totally split over $\PP^1_{\hat\ell}$.  So the fibre of $Y \to \PP^1_{\hat\ell}$ over any 
$\hat P'_i$ or $\hat P''_j$ consists of $|H|$ distinct $\hat\ell$-points, where $H$ is the Galois group of $Y$ over $\PP^1_{\hat\ell}$.
Note that the points in $B$ have the property that the closures of their images under $\bar Y \to \PP^1_R$ each meet the closed fibre $\PP^1_k$ precisely at the point $P_0$.  Similarly, the points of $B_0$ have the property that the closures of their images under $\bar Y \to \PP^1_S$ each meet the closed fibre of $\PP^1_S$ precisely at the point $P'_0$, over the closed point $P_0$ of $\PP^1_R$.  So 
we may index the points in the fibres of $Y \to \PP^1_{\hat\ell}$ over 
$\hat P'_i$ and $\hat P''_j$ as $\hat Q'_{i,\delta}$ and $\hat Q''_{j,\delta}$ for $1 \le \delta \le |H|$, such that the closures of $\hat Q'_{i,1}$ and $\hat Q''_{j,1}$ in $\bar Y$ pass through $Q$.

Let $\Sigma$ be the reduced intersection of the closure of $B$ in $\bar Y$ with the closed fibre of $\bar Y$.
For each non-negative integer $\iota$, consider the $\iota^{\rm th}$ blow-up $(\bar Y_\iota, B, \Sigma_\iota)$ of $(\bar Y, B, \Sigma)$.  
Lemma~3.2 applies to $\bar Y$ and Lemma~3.3 applies to the $E$-Galois cover $\bar Y \to \PP^1_R$; so the conclusions of these lemmas hold for sufficiently large $\iota>h$.  Fix such a value $\iota = \lambda$, write $\tilde Y$ for $\bar Y_\lambda$, and write $\Delta'$ for $\Delta_\lambda$.  Thus the action of $E$ on $\bar Y$ lifts to an action of $E$ on $\tilde Y$.  So the Galois groups $G \subset E$ of $\bar Y \to \bar X$ and $H \subset E$ of $\bar Y \to \PP^1_S$ also lift to actions on $\tilde Y$.  

Let $\tilde X = \tilde Y/G$, $\tilde \PP_R = \tilde Y/E$, and $\tilde \PP_S = \tilde Y/H$.  Thus $\tilde Y \to \tilde \PP_S$ is an $H$-Galois cover, and there is a birational proper morphism $\omega:\tilde \PP_S \to \PP^1_S$ such that the composition $\tilde Y \to \bar Y \to \PP^1_S$ factors as $\tilde Y \to \tilde \PP_S \to \PP^1_S$.  
Similarly, the composition $\tilde Y \to \bar Y \to \PP^1_R$ factors as $\tilde Y \to \tilde \PP_R \to \PP^1_R$, where $\tilde \PP_R \to \PP^1_R$ is a birational proper morphism.
Also, the general fibre of $\tilde \PP_S$ is $\PP^1_{\hat\ell}$
and that of $\tilde \PP_R$ is $\PP^1_{\hat k}$.  Let $P_i', P_j''$ be the closed points of $\tilde \PP_S$ where the closures $\bar P'_i, \bar P''_j$ of $\hat P'_i, \hat P''_j$ respectively meet the closed fibre of $\tilde \PP_S$.  Then these $r+s$ closed points, and their $k$-conjugates, are all distinct, because the fibres of $\tilde Y \to \tilde \PP_S$ over $\hat P'_i, \hat P''_j$ and their $\hat k$-conjugates have disjoint closures.  Using Lemma~3.2(iii), observe that the closed fibre $(t=0)$ of $\tilde \PP_S$ is locally irreducible at each of the points $P_i',P_j''$, since the corresponding assertion holds for the points of $\tilde Y$ over $P_i',P_j''$ (these being points of $\Delta'$).

\smallskip

\noindent{\it Step 3.} Construction of the points $\hat P_i'^*,\hat P_j''^*$:

\smallskip

Since the points of $Y$ over a given $\hat P_i'$ (which form a subset of $B_0 \subset B$) have disjoint closures in $\tilde Y$, it follows that the fibre of $\tilde Y \to \tilde \PP_S$ over 
$P_i'$ consists of $|H|$ distinct closed points; we may label these as $Q_{i,\delta}'$, for $1 \le \delta \le |H|$, 
where $Q_{i,\delta}'$ is the point of $\tilde Y$ over $P_i'$ that is in the closure of $\hat Q_{i,\delta}' \in Y$ in $\tilde Y$.  Similarly,
the fibre of $\tilde Y \to \tilde \PP_S$ over each $P_j''$
consists of $|H|$ distinct closed points; and we may index these as
$Q_{j,\delta}''$, for $1 \le \delta \le |H|$, such that the closure of $\hat Q_{j,\delta}''$ in $\tilde Y$ meets the closed fibre at $Q_{j,\delta}''$.  Since the $\hat\ell$-points $\hat Q_{i,\delta}', \hat Q_{j,\delta}''$
and their $\hat k$-conjugates have disjoint closures in $\tilde Y$, it follows that the closed points 
$e(Q_{i,\delta}')$ and $e(Q_{j,\delta}'')$ are all distinct, where $1 \le i \le r$, $1 \le j \le s$, $1 \le \delta \le |H|$, and where $e \in E = {\rm Gal}(Y/\PP^1_{\hat k})$ ranges over a set of coset representatives of $E/H$.

If $h'>\lambda$ is sufficiently large, then any $\hat \ell$-point 
$\hat P_i'^*$ [resp.\ $\hat P_j''^*$]
of $\PP^1_{\hat \ell}$ that is congruent to $\hat P_i'$ [resp.\ $\hat P_j''$] modulo $t^{h'}$ (relative to $\PP^1_S$) will have the property that its closure in 
$\tilde \PP_S$ will contain $P_i'$ [resp.\ $P_j''$], and the two Cartier divisors $\omega^*(\bar P_i')$ and $\omega^*(\bar P_i'^*)$ [resp.\ $\omega^*(\bar P_j'')$ and $\omega^*(\bar P_j''^*)$] have the same multiplicities at each generic point of the closed fibre of $\tilde \PP_S$.  (Here $\bar P_i'^*$ and $\bar P_j''^*$ are the closures of $\hat P_i'^*$ and $\hat P_j''^*$ in $\PP^1_S$.)  Applying Lemma~3.4(b) to the points $P_i'$ and $P_j''$ for $1 \le i \le r$ and $1 \le j \le s$ (with $h'$ playing the role of the $h$ in the lemma), we obtain $\hat \ell$-points $\hat P_i'^*, \hat P_j''^*$ congruent to $\hat P_i', \hat P_j''$ modulo $t^{h'}$ and
satisfying the conclusions there.  
Since there are 
{\new ${\rm card}\,\hat k$}
choices for these points, 
they may be chosen so that they and their $\hat k$-conjugates are distinct from the points $\hat P_i', \hat P_j''$ and their $\hat k$-conjugates{\new, and in (e) so that $P_0' \times_\ell \hat \ell$ is not among these points.
Since there are ${\rm card}\,\hat k$ ways of choosing the finite set of points $\hat P_i', \hat P_j'', \hat P_i^*$ satisfying the above conditions, there is a set of ${\rm card}\,\hat k$ disjoint choices of that set.} 

\smallskip

\noindent{\it Step 4.} Construction of the rational functions $\varphi_i, \psi_j$:

\smallskip

According to the conclusion of Lemma~3.4(b), $\hat P_i'^*$ totally splits in the Galois cover $Y \to \PP^1_{\hat\ell}$; so each point of $\tilde Y$ over $P_i'$ is in the closure of a unique point of $Y$ over $\hat P_i'^*$.
Thus we may index the points of $Y$ over $\hat P_i'^*$ as 
$\hat Q_{i,\delta}'^*$ for $1 \le \delta \le |H|$, such that the closure of $\hat Q_{i,\delta}'^*$ in $\tilde Y$ meets the closed fibre at $Q_{i,\delta}'$ (and thus meets the closure of $\hat Q_{i,\delta}'$
there).  The corresponding assertions hold for $\hat P_j''^*$, $\hat P_j''$, $\hat Q_{j,\delta}''^*$, $\hat Q_{j,\delta}''$.

For $1 \le i \le r$, the two $S$-points $\bar P_i', \bar P_i'^*$ are linearly equivalent on $\PP^1_S$.  So there is a rational function $\varphi_i$ on $\PP^1_S$ whose divisor is $\bar P_i' - \bar P_i'^*$.  
Viewing $\varphi_i$ as a rational function on $\tilde \PP_S$ via $\omega:\tilde \PP_S \to \PP^1_S$, its divisor is $\omega^*(\bar P_i') - \omega^*(\bar P_i'^*)$; so its zero and pole divisors on $\tilde \PP_S$ are the closures of $\hat P_i'$ and $\hat P_i'^*$ in $\tilde \PP_S$ (since the supports 
of $\omega^*(\bar P_i')$ and $\omega^*(\bar P_i'^*)$
on the closed fibre are equal and hence cancel).
Hence as a rational function on $\tilde \PP_S$, $\varphi_i$ is defined and invertible on the closed fibre except possibly at the one point $P_i'$.  But the closed fibre is a connected projective curve.  So actually $\varphi_i$ defines a non-zero constant function on this closed fibre; and after multiplying it by a unit in $S$, we may assume that this constant is $1$.  Thus $\varphi_i$ restricts to the constant function $1$ at each generic point of the closed fibre of $\tilde \PP_S$.

Similarly, for $1 \le j \le s$ there is a rational function $\psi_j$ on $\tilde \PP_S$ restricting to the constant function $1$ at each generic point of the closed fibre, with zero and pole divisors being the closures of $\hat P_j''$ and $\hat P_j''^*$ respectively.  

\smallskip

\noindent{\it Step 5.} Construction the rational functions $\Phi_i, \Psi_j$ and the cover $Z \to Y \to X$:

\smallskip

Taking norms, let $\Phi_i = N_{\hat \ell/\hat k}\,\varphi_i$; this is the product of $\varphi_i$ and its $\hat k$-conjugates.  Thus $\Phi_i$ defines a rational function on $\PP^1_R$ and hence on $\tilde \PP_R$, where its zero divisor is the closure of $\hat P_i'$ and its conjugates, and where its pole divisor is the closure of $\hat P_i'^*$ and its conjugates.  Moreover $\Phi_i$ restricts to the constant function $1$ at each generic point of the closed fibre (since this holds for its pullback to $\tilde \PP_S$).
Similarly, for $1 \le j \le s$ we let $\Psi_j$ be the rational function on $\tilde \PP_R$ obtained by taking the norm of $\psi_j$.  This has the corresponding properties, with respect to the points $\hat P_j'',\hat P_j''^*$.

The functions $\Phi_i, \Psi_j$ on $\tilde \PP_R$ pull back to rational functions on $\tilde X$, which we again denote by $\Phi_i$ and $\Psi_j$.  These functions again restrict to the constant function $1$ at the generic points of the closed fibre. Their zero and pole divisors are the inverse images of those on $\tilde \PP_R$, and are therefore reduced, since the supports of those divisors do not lie in the branch locus of $X \to \PP^1_{\hat \ell}$.  

So Proposition~3.5 applies, yielding an irreducible normal $N$-Galois cover $\tilde Z \to \tilde Y$ such that $\tilde Z \to \tilde X$ is $\Gamma$-Galois{\new, and satisfying conditons (i)-(vi) of that result}.  The generic fibre $Z \to Y$ of this cover is branched precisely at the supports on $Y$ of the divisors of the $\Phi_i$'s and the zero divisors of the $\Psi_j$'s; i.e.\ at the points of $Y$ over $\hat P'_i, \hat P_i'^*, \hat P''_j \in \PP^1_{\hat \ell}$.  Moreover $\hat \ell$ is the algebraic closure of $\hat k$ in the function field of $Z$.  So $Z \to Y$ provides a proper regular solution to the given embedding problem, proving (a) of the theorem in the general case.

\smallskip

\noindent{\it Step 6.}  Verification of properties (b)-{\new (e).}

\smallskip

The cover $\tilde Z \to \tilde Y$ provided by Proposition~3.5 has the property that the points of $D \subset Y \subset \tilde Y$ are totally split, and that the decomposition groups at the points of $Z$ over $\delta\in D$ are the conjugates of $\sigma(G_\delta)$, where $\sigma$ is the given splitting of the exact sequence and $G_\delta$ is the decomposition group of $Y \to X$ at $\delta$.  
{\new For the last part of (b), consider an \'etale subcover 
$W \to Y$ of $Z \to Y$.  The Galois group $\Gal(Z/W)$, which is a subgroup of $N = \Gal(Z/Y)$, must contain all the inertia groups of $Z \to Y$.  But by condition (iv) of Proposition~3.5, each $n_i$ [resp.\ $m_j$] generates an inertia group of $Z_0 \to Y$ over $\hat P'_i$ [resp.\ $\hat P''_j$].  Since the $r+s$ elements $n_i$ and $m_j$ generate $N$, if follows that $\Gal(Z/W)$ is all of $N$; i.e., $W = Y$.}  So (b) of the theorem holds.

For {\new the first part of} (c), we want to show that the number of proper regular solutions constructed above to the given embedding problem ${\cal E} = (\alpha:G_K \to G, f:\Gamma \to G)$ is equal to the cardinality of $\hat{k}=k((t))$.  As before, it suffices to show that the cardinality of the set of proper regular solutions to $\cal E$ is at least ${\rm card}\, \hat k$.

The above construction of a solution $Z \to Y$ depended on a choice of $2r+s$ points, $\hat P'_i, \hat P_i'^*, \hat P''_j$ (where $1 \le i \le r$ and $1 \le j \le s$).  By Lemma~3.4(b), there are ${\rm card}\,\hat k$ choices for these points (using also that $r$ and $s$ are not both $0$ in part (c)).  Distinct choices for this set of points will yield non-isomorphic solutions $Z \to Y$, since the branch locus of $Z \to Y$ is precisely the inverse image of these points $\hat P'_i, \hat P_i'^*, \hat P''_j \in \PP^1_{\hat \ell}$ under $Y \to \PP^1_{\hat \ell}$ (using also that the chosen generators of $N$ are all non-trivial).  So there are at least ${\rm card}\,\hat k$ choices for the cover $Z \to Y$, proving {\new the first part of} (c).

{\new For the second part of (c), first note that there exists a set of $m$ choices $\hat P'_{i,\alpha}, \hat P_{i,\alpha}'^*, \hat P''_{j,\alpha}$ of the points $\hat P'_i, \hat P_i'^*, \hat P''_j$ such that no two of these points (as $\alpha$ varies over the index set $S$ of cardinality $m$) are equal.  For each $\alpha \in S$, choose as above a corresponding proper solution $Z_\alpha \to Y$.  It suffices to show that the solutions $Z_\alpha \to Y$ are {\new independent}.  By induction, this would follow from showing that if $\alpha_0,\dots,\alpha_n \in S$, then the cover $Z_0 := Z_{\alpha_0} \to Y$ is linearly disjoint from the minimal cover $Z_1 \to Y$ that dominates the covers $Z_{\alpha_1},\dots, Z_{\alpha_n}\to Y$.  
Since these covers are Galois, this is equivalent to showing that $Z_0, Z_1 \to Y$ dominate no common non-trivial intermediate cover of $Y$.  So suppose that $W \to Y$ is a connected normal cover that is dominated by $Z_0 \to Y$ and $Z_1 \to Y$.  By the choice of the covers $Z_\alpha \to Y$, the branch loci of $Z_0 \to Y$ and $Z_1 \to Y$ are disjoint; hence $W \to Y$ is \'etale.  By the last assertion of part (b), $W=Y$.  This shows that $Z_0$ and $Z_1$ are linearly disjoint over $Y$.}

To verify the conclusion of part (d), {\new it suffices to show that the ${\rm card}\,\hat k$ {\new independent} proper solutions constructed above have the asserted properties.  So consider one of these proper solutions $Z \to Y$, and let $\tilde Z \to \tilde Y$ be as in the above construction.  Also,}
let $\bar Z$ be the normalization of $\bar Y$ in $Z$.  So the composition $\tilde Z \to \tilde Y \to \bar Y$ factors as $\tilde Z \to \bar Z \to \bar Y$, where $\tilde Z \to \bar Z$ is a proper birational morphism (contracting components of the closed fibre).   Since the closure of $\hat Q'_{i,1}$ [resp.\ $\hat Q''_{j,1}$] in $\tilde Y$ meets the closed fibre at $Q'_{i,1}$ [resp.\ $Q''_{j,1}$], and since the corresponding closure in 
$\bar Y$ meets the closed fibre just at $Q$, it follows that each $Q'_{i,1}$ and each $Q''_{j,1}$ must map to $Q$ under $\tilde Y \to \bar Y$.  Let $\bar Y_0$, $\tilde Y_0$ be the closed fibres of $\bar Y$, $\tilde Y$.
Since the closed fibre of $\tilde Z \to \tilde Y$ is a mock cover by Proposition~3.5, the same is true for the closed fibre of $\bar Z \to \bar Y$, with the identity copy of $\tilde Y_0$ in the closed fibre of $\tilde Z$ mapping to the identity copy of $\bar Y_0$ in the closed fibre of $\bar Z$.  So the points $Q'_{i,1}$ and $Q''_{j,1}$ on the identity copy of $\tilde Y_0$ must map to the point 
$Q$ on the identity copy of $\bar Y_0$. 
But the inertia groups of $\tilde Z \to \tilde Y$ at the points $Q'_{i,1}$ and $Q''_{j,1}$ on the identity sheet are respectively generated by $n_i$ and $m_j$.  So the inertia group of $\bar Z \to \bar Y$ at the point $Q$ on the identity sheet contains every $n_i$ and $m_j$.  But these elements generate $N = \Gal(\bar Z/\bar Y)$.  So $\bar Z \to \bar Y$ is totally ramified over $Q$.  Thus each of the ${\rm card}\,\hat k$ {\new independent} {\new proper} solutions constructed above have this property, in addition to (a) and (b).  So the first part of (d) of the theorem holds. 

For the local part of (d), consider a global solution $Z \to Y$, and its pullback $Z^* \to Y^* = {\rm Spec}\,{\hat {\cal O}}_{{\bar Y},Q} \to {\bar Y}$.   The cover $Z \to Y$ is branched non-trivially at the points of $Y$ lying over $\hat P'_i, \hat P_i'^*, \hat P''_j \in \PP^1_{\hat \ell}$; in particular, at the $\hat \ell$-points $\hat Q'_{i,1}, \hat Q'^*_{i,1}, \hat Q''_{j,1} \in Y$.  But the closures of those latter points in $\bar Y$ contain the closed point $Q \in \bar Y$.   Hence the branch locus of $Z^* \to Y^*$ contains the closures of the (distinct) $\hat \ell$-points $\hat Q'_{i,1}, \hat Q'^*_{i,1}, \hat Q''_{j,1} \in Y^*$.  Since there are ${\rm card}\,\hat k$ choices for $\hat P'_i, \hat P''_j$, the same is true for $\hat Q'_{i,1}, \hat Q'^*_{i,1}, \hat Q''_{j,1}$ and hence for $Z^* \to Y^*$.  So there are at least (and hence exactly) ${\rm card}\, \hat k$ such pullbacks.  
{\new Moreover these pullbacks are linearly disjoint by the same argument as in the proof of the second part of (c).}  Each of them is connected because $\bar Z$ and hence $Z^*$ has only one closed point over $Q$, and each connected component of $Z^*$ must contain a point over $Q$ (as $Z^* \to Y^*$ is finite).  Since $\bar Z$ and hence $Z^*$ is normal, it follows that $Z^*$ is irreducible.  

Similarly, if $Q$ is the only point of $\bar Y$ over $P \in \bar X$, then there is just one point of $\bar Z$ over $P$, and so the pullbacks of $\bar Y  \to \bar X$ and of $\bar Z \to \bar X$
under $X^* = {\rm Spec}\,{\hat {\cal O}}_{{\bar X},P} \to {\bar X}$ are also irreducible.  So in that case these pullbacks are the same as the above $Y^*$ and $Z^*$.  Since distinct covers $Z^* \to Y^*$ remain distinct upon composition with $Y^* \to X^*$, it follows that there are exactly ${\rm card}\, \hat k$ distinct covers $Z^* \to X^*$ arising from these pullbacks.  This proves the local part of (d).  {\new

Finally, for assertion (e), recall that $Z \to Y$ is ramified precisely at the points of $Y$ over $\hat P'_i, \hat P_i'^*, \hat P''_j \in \PP^1_{\hat \ell}$.  These points were chosen so as not to include $P_0' \times_\ell \hat \ell$.  Also, the closures of these points of $\PP^1_\ell$ each meet the closed fibre of $\PP^1_S$ precisely at $P_0'$.  So assertion (e) follows.}
\qed

\medskip
 
Next, we extend  Theorem~4.1 to arbitrary large fields, not necessarily  of the form $k((t))$.  Recall from [Po96] that a field $F$ is called {\it large} if it has the property that every smooth $F$ curve with an $F$-point has infinitely many $F$-points.  According to [Po96, Proposition 1.1] this is equivalent to each of the following two conditions:  For every smooth integral $F$-variety with an $F$-point, the set of such points is dense.  Every $F$-variety with an $F((t))$-point has an $F$-point.  

Large fields are infinite, and examples of large fields include algebraically closed fields and fields of the form $\hat k = k((t))$.  (The case of $\hat k$ follows by choosing a covering map from a given $\hat k$-curve to the line, such that the given $\hat k$-point is unramified; and then using Hensel's Lemma to choose infinitely many $t$-adically nearby $\hat k$-points.)  In fact these examples satisfy a stronger condition:

We will say that a field $F$ is {\it very large} if it has the property that for every smooth $F$-curve with an $F$-point, the set of $F$-points has cardinality equal to ${\rm card}\,F$.   Trivially, every very large field is large.  Also, algebraically closed fields and fields of the form $\hat k = k((t))$ are very large (by the same arguments as for large).

We then have the following generalization of parts (a)-(c) of Theorem~4.1, which also extends the results of Pop [Po96, Main Theorem~A] and Haran-Jarden [HJ, Theorem~6.4].  Those results assert the existence of proper regular solutions for embedding problems for curves over large fields, and assert that such solutions can be chosen to be totally split over one given unramified point.  (See also [Ha03, Remark~5.1.11].)
  
\medskip

\noindent{\bf Theorem~4.3.} {\sl Let $F$ be a large field and let $X$ be a smooth projective connected $F$-curve, with function field $K$.  Let ${\cal E} = (\alpha:G_K \to G, f:\Gamma \to G)$ be a non-trivial finite split embedding problem with section $\sigma$, let $\pi:Y \to X$ be the $G$-Galois branched cover of normal curves corresponding to $\alpha$, and let $D \subset Y$ be a finite set of closed points.  

a) Then there are infinitely many distinct proper regular solutions to $\cal E$ such that the corresponding cover 
$Z \to Y$ is totally split over the points of $D${\new ;}
the decomposition groups of $Z \to X$ at the points of $Z$ over $\delta \in D$ are the conjugates of $\sigma(G_\delta)$, where $G_\delta$ is the decomposition group of $Y \to X$ at $\delta${\new; and $Z \to Y$ has no non-trivial \'etale subcover $W \to Y$}.

b) Moreover, if $F$ is very large, then the cardinality of the set of such solutions is equal to ${\rm card}\,F${\new, and there is  {\new an independent} set of such solutions having cardinality equal to ${\rm card}\,F$}.}

\medskip

\noindent{\it Proof.}    
Taking a separating transcendence basis for $K$ over $F$, we obtain a cover $X \to \PP^1_F$.  The field $F$ is infinite, since it is large; so there is an $F$-point of $\PP^1_F$ that is not in the branch locus of the composition $Y \to X \to \PP^1_F$.  Let $Q$ be a point of $Y$ lying over this point of $\PP^1_F$, and consider its image in $X$.  So the residue field at this image is separable over $F$.  

Let $R = F[[t]]$ and $\hat F = F((t))$, and let $\hat X = X \times_F \hat F$ and $\hat Y = Y \times_F \hat F$.  So the function field of $\hat X$ is $\hat K := K \otimes_F \hat F$.  The cover $Y \to X$ induces a
$G$-Galois connected normal cover $\hat Y \to \hat X$, and this in turn corresponds to an epimorphism $\hat\alpha:G_{\hat K} \to G$.  We thus obtain a split embedding problem $\hat \E = (\hat\alpha:G_{\hat K} \to G, f:\Gamma \to G)$ for $\hat K$.
For every closed point $P$ of $Y$, we may consider the closed point $\hat P = P \times_F \hat F$ of $\hat Y$.  Let $\hat D \subset \hat Y$ be the finite subset consisting of the points $\hat P$ for each $P \in D$ together with the point $\hat Q$.  Let $\bar X = X \times_F R$ and let $\bar Y = Y  \times_F R$.  Since $X$ is smooth over $F$, it follows that $\bar X$ is smooth over $R$.  Similarly, $Y$ is a regular scheme (being a normal curve), and hence so is $\bar Y$.
Regarding $Y$ as the closed fibre of $\bar Y$, we may view $Q$ as a point of that closed fibre.  

We may now apply Theorem~4.1 to the $\hat F$-curve $\hat X$, the finite split embedding problem $\hat \E$, the set $\hat D \subset \hat Y$, the $R$-model $\bar Y \to \bar X$, and the closed point $Q \in \bar Y$.   

By Theorem~4.1, there is a proper regular solution $\hat\beta:G_{\hat K} \to \Gamma$ such that the corresponding $\Gamma$-Galois cover $\hat Z \to \hat X$ dominates $\hat Y \to \hat X$ and satisfies conditions (b) and (d) there (with respect to $\hat F$, $\hat D$, $\bar Y \to \bar X$, and $Q \in \bar Y$).   Let $\bar Z$ be the normalization of $\bar X$ in $\hat Z$.  

The scheme $\bar Y$ is regular and the scheme $\bar Z$ is normal; so Purity of Branch Locus applies to $\bar Z \to \bar Y$. 
{\new By Theorem~4.1(e), no irreducible component of the branch locus of $\bar Z \to \bar Y$ lies over a divisor on $\bar X$ that is induced from a divisor on its closed fibre by base change from $F$ to $R$.  Thus no point in the branch locus of $\hat Z \to \hat Y$ is induced from an $F$-point of $Y$ by base change to $\hat F$.}
Let $r$ be the degree of this branch locus. 

The cover $\hat Z \to \hat Y$ descends from $\hat F$ to some finite type $F$-subalgebra $A$ of $\hat F$.  Thus for some such $A$, if we let $X_A = X \times_F A$ and $Y_A = Y \times_F A$, there is a connected projective normal $A$-curve $Z_A$ and a covering morphism $Z_A \to Y_A$ such that the composition $Z_A \to X_A$ is $\Gamma$-Galois; the cover $Z_A \to Y_A$ is totally split over every point of $D_A := D \times_F A$, with $\sigma(G_\delta)$ a decomposition group of $Z_A \to X_A$ over any point of $\delta \times_F A$ for $\delta \in D$; {\new $Z_A \to Y_A$ has no non-trivial \'etale subcover $W_A \to Y_A$;} and the pullback $Z_A \times_A \hat F \to \hat Y$ is isomorphic to $\hat Z \to \hat Y$.  Necessarily, {\new no component of}
the branch locus of $Z_A \to Y_A$ is induced from $F$, because of this property for $\hat Z \to \hat Y$.  That is, the branch loci of the fibres of $Z_A \to Y_A$ over the points of $\Spec A$ are non-constant {\new and vary without base points}.  After replacing $A$ by a basic affine open subset, we may assume that $\Spec A$ is smooth over $F$.  Also $\Spec A$ is irreducible, since $A \subset \hat F$ is a domain.

The curve $\hat Z$ is geometrically irreducible as a scheme over the algebraic closure of $\hat F$ in the function field of $\hat Y$, since $\hat Z \to \hat X$ corresponds to a proper regular solution to the embedding problem $\hat \E$.  So the general fibre of $Z_A$ is geometrically irreducible as a scheme over the algebraic closure of ${\rm frac}(A)$ in the function field of $Y_A$.
But the set of points of $\Spec A$ at which the fibre of $Z_A$ is geometrically irreducible is a constructible set [Gr66, 9.7.7]{\new, as is the set of points at which the fibre is geometrically normal [Gr66, 9.9.4]}; so after replacing $A$ by a basic affine open subset, 
we may assume that every closed fibre of $Z_A \to \Spec A$ is geometrically irreducible {\new and normal}.  We may also assume that each such fibre is generically \'etale {\new over the corresponding fibre of $Y_A$; and that} 
each fibre has exactly $r$ geometric branch points{\new;} 
each {\new of which is} branched non-trivially{\new, with the same inertia groups as at the corresponding geometric branch points of $Z \to Y$.  Since $Z \to Y$ has no non-trivial subcovers $W \to Y$ by part (b) of Theorem~4.1, its inertia groups generate $N = \Gal(Z/Y)$; and hence the same is true for each fibre of $Z_A \to Y_A$.}
So for every $F$-point of $\Spec A$, 
the fibre of $Z_A \to Y_A$ over this point will be a proper regular solution to the given embedding problem with desired properties.    

It remains to show that the set of such solutions has the asserted cardinality in parts (a) and (b) above.

Let $Y^{(r)}$ be the $r^{\rm th}$ symmetric power of $Y$, and let $b:\Spec A \to Y^{(r)}$ be the morphism assigning to each point of $\Spec A$ the branch locus of the corresponding normalized geometric fibre of $Z_A \to Y_A$.  Thus $b$ is non-constant, and so its image $B$ has dimension $\ge 1$.  Let $S$ be the set of $F$-points of $\Spec A$.  {\new Since $F$ is large and since $\Spec A$ is smooth and contains an $\hat F$-point (viz.\ the point corresponding to the inclusion $A \hookrightarrow \hat F$), it contains an $F$-point, and in fact the set $S$ is dense in $\Spec A$.}  Hence $b(S)$ is dense in $B$, and is thus infinite.  Since covers with distinct branch loci are non-isomorphic, it follows that the set of proper regular solutions as above is infinite.  This proves part (a).

For part (b), we wish to show in the very large case that the set of such solutions has cardinality equal to ${\rm card}\,F${\new, and that there is {\new an independent}  set of such solutions that has cardinality ${\rm card}\,F$}.  As before, {\new for the first of these assertions} it suffices to show that there are at least that many {\new such solutions; and this would follow from the second assertion}.  Proceeding as above, take an $F$-point $P$ on $\Spec A$ (of which there are infinitely many since $F$ is large).  
{\new Let $B$ be the subset of $Y^{(r)}$ consisting of the effective divisors on $Y$ whose support is not disjoint from that of $b(P)$.
Then $b^{-1}(B)$} 
is closed and is strictly contained in $\Spec A$, because 
{\new the branch loci of the fibres of $Z_A \to Y_A$ forms a base-point free family}.  
Since $\Spec A$ is smooth and connected, there is a smooth curve $C$ on $\Spec A$ that passes through $P$ and is not contained in 
{\new $b^{-1}(B))$.} 
{\new Hence for each point $Q'$ of $Y$, there are at most finitely many points of $C$ for which the $Q'$ is a branch point of the corresponding fibre of $Z_A \to Y_A$.}
Since $F$ is very large, the set of $F$-points on $C$ has cardinality ${\rm card}\,F$.  Hence the set of proper regular solutions induced by the $F$-points of $C$ has cardinality equal to ${\rm card}\,F$; 
{\new and there exists a set $S$ of $F$-points of $C$ having cardinality ${\rm card}\,F$ and whose fibres of $Z_A \to Y_A$ have pairwise disjoint branch loci.} 
{\new But as shown above, each fibre of $Z_A \to Y_A$ has the property that its inertia groups generate its Galois group $N$; and thus it has no non-trivial \'etale subcovers.  Hence the covers corresponding to the points of $S$ are linearly disjoint, by the same argument used to prove the second part of condition 4.1(c) in the proof of Proposition~4.2.}
This proves part (b).
\qed

\medskip

\noindent {\bf Corollary~4.4.} {\sl If $K$ is the function field of a smooth projective curve over a very large field $k$, then the absolute Galois group of $K$ is {\new semi-free}}.

\medskip

\noindent{\it Remark.} If $k$ is also perfect, then every field extension $K$ of transcendence degree $1$ over $k$ is the function field of a smooth projective $k$-curve, and hence Corollary~4.4 applies.

\medskip

Because of Corollary~4.4, it {\new is natural to ask}

\smallskip

\noindent{\bf Question~4.5.} Is every large field very large?  

{\new \smallskip 

The answer to this question is now known to be yes, by a result of F.~Pop; see Proposition~3.3 of [Ha08] or {\new Theorem~5 of [HS06]}.  Therefore, by Corollary~4.4, the absolute Galois group of the function field of a smooth projective curve over a large field is {\new semi-free}.}

\bigskip

\noindent {\bf Section 5. Main Result}

\medskip

In this section we prove the main result of the paper, stated as Theorem~1.1 in the introduction.

\medskip

\noindent {\bf Theorem~5.1.}  (Main theorem) {\sl Let $k$ be a field. 
Then the absolute Galois group of the field $K := k((x,t))$ is 
{\new semi-free} of rank ${\rm card}\, K$.  Equivalently, 
every non-trivial finite split embedding problem for $K$ has {\new exactly} ${\rm card}\, K$ {\new non-isomorphic} proper solutions{\new , including a set of ${\rm card}\, K$ {\new independent}  proper solutions.}}

\medskip

This result generalizes Theorem~5.3.9 of [Ha03], which 
showed the existence of proper solutions in 
the case $k={\C}$, but did not consider the number of such solutions.  
The proof of that earlier result relied heavily on the fact that $\C$ is algebraically closed and of characteristic $0$, in order to know that every finite extension of $\C((x))$ is obtained by adjoining some $n^{\rm th}$ root of $x$, and in order to be able to use Abhyankar's Lemma.  Here, under the more general hypothesis that $k$ is an arbitrary field, the proof of Theorem~5.1 is completely different, and relies on Theorem~4.1.  

Theorem~5.1 follows immediately from Proposition~5.3 below, a geometric version of the main theorem stated in terms of covers.
Proposition~5.3 considers a $G$-Galois branched cover $Y^* \to X^*=\Spec k[[x,t]]$ and
a split short exact sequence $1 \to N \to \Gamma \ {\buildrel f \over \to} \ G \to 1$.
This proposition asserts that there exist ``many'' $\Gamma$-Galois covers $Z^* \to X^*$ dominating $Y^* \to X^*$ and having an additional splitting property.  
The proof of Proposition~5.3 uses Theorem~4.1 to solve an embedding problem for the $k((t))$-line $X$; then we close that up over the $k[[t]]$-line $\bar X$ (i.e.\ normalizing $\bar X$ in these covers of $X$); and finally we restrict from $\bar X$ to $X^*$, which we regard as the complete local neighborhood of a closed point of $\bar X$.  In the process, we rely on a preliminary result, Lemma~5.2, which permits us to pass from an embedding problem over $X^*$ to a more global embedding problem in Proposition~5.3.

\medskip

In the following result, we consider the morphism from $X^* = \Spec k[[x,t]]$ to the affine $x$-line $\AA^1_{k[[t]]}$ that corresponds to 
the inclusion of rings $k[[t]][x] \hookrightarrow k[[x,t]]$.  Composing with the inclusion $\AA^1_{k[[t]]} \hookrightarrow \bar X := \PP^1_{k[[t]]}$, we obtain a morphism $X^* \to \bar X$.

\medskip

\noindent{\bf Lemma~5.2.} {\sl Let $k$ be a field of characteristic $p \ge 0$, let $G$ be a finite
group, and let $Y^* \to X^* =
\Spec k[[x,t]]$ be a connected normal $G$-Galois cover that
is unramified over the generic point of $(t=0)$.
Then there is a connected normal generically smooth $G$-Galois cover $\bar Y \to \bar X$ of $k[[t]]$-curves whose pullback via $X^* \to \bar X$ is the given cover $Y^* \to X^*$.}

\medskip

\noindent{\it Proof.}  Let $Y^*_0 \to X^*_0 = \Spec k[[x]]$ be the reduction modulo $t$ of $Y^* \to X^*$.  By hypothesis, this finite morphism is generically \'etale.
Let $Y^\circ_0 \to X^\circ_0  = \Spec k((x))$ be the generic fibre of $Y^*_0 \to X^*_0$, and let $Z^\circ_0$ be a connected component of $Y^\circ_0$ (corresponding to a finite field extension of $k((x))$).  So $Z^\circ_0 \to X^\circ_0$ is $H$-Galois for some $H \subset G$, and $Y^\circ_0 = \Ind_H^G\, Z^\circ_0$.  Applying the theorem of Katz-Gabber [Ka, Theorem 1.4.1], there is an \'etale $H$-Galois cover $Z''_0 \to X''_0 := \PP^1_k - \{0,\infty\} = \Spec k[x,x^{-1}]$ whose pullback to $X^\circ_0$ is $Z^\circ_0 \to X^\circ_0$.  Let $Z'_0 \to X'_0 := \PP^1_k - \{0\} = \Spec k[x^{-1}]$ be the normalization of $X'_0$ in $Z''_0$.  Thus the pullback of $Z'_0 \to X'_0$ to $X^\circ_0$ is again $Z^\circ_0 \to X^\circ_0$; and so we may identify the pullback of $Y'_0 := \Ind_H^G\, Z'_0 \to X'_0$ to $X^\circ_0$ with $Y^\circ_0 \to X^\circ_0$ as a (disconnected) $G$-Galois cover.  

Meanwhile, let $\bar X^\circ = \Spec k((x))[[t]]$ and let $\bar Y^\circ = Y^* \times_{X^*} \bar X^\circ$.
Since $Y^* \to X^*$ is unramified at the 
generic point of $(t=0)$, its pullback $\bar Y^\circ \to \bar X^\circ$ 
under $\bar X^\circ \to X^*$
is unramified at the closed point $(t=0)$ of $\bar X^\circ$ and is therefore \'etale.  So by Hensel's Lemma, we have that $\bar Y^\circ = Y^\circ_0 \times_{X^\circ_0} \bar X^\circ$ canonically, as $G$-Galois covers (compatibly with the restriction of $\bar Y^\circ$ to $Y^\circ_0$).  

Let $\bar X' = \Spec k[x^{-1}][[t]]$ and let $\bar Y' = Y'_0 \times_{X'_0} \bar X'$.  So we have identifications of the $G$-Galois covers $\bar Y' \times_{\bar X'} \bar X^\circ = Y'_0 \times_{X'_0}  \bar X^\circ = Y^\circ_0 \times_{X^\circ_0} \bar X^\circ = \bar Y^\circ$.  
Thus the pullbacks of $Y^* \to X^*$ and of $\bar Y' \to \bar X'$ to $\bar X^\circ$ are each identified with $\bar Y^\circ \to \bar X^\circ$.

So by formal patching ([HS99, \S1, Cor.\ to Thm.~1], [Pr, Thm.~3.4]), there is a unique $G$-Galois cover $\bar Y \to \bar X$ whose pullbacks to $X^*, \bar X', \bar X^\circ$ respectively agree with $Y^* \to X^*$, $\bar Y' \to \bar X'$, and $\bar Y^\circ \to \bar X^\circ$, compatibly with the above identifications.  Now $Y^*$ is connected; $Y^*$ and $\bar Y'$ are normal; and $Y'_0 = \bar Y \times_{\bar X} X'_0$ is generically smooth over $k$.  So $\bar Y$ is connected and normal, and is generically smooth over $k[[t]]$.  
\qed

\bigskip

We now prove the geometric form of our main result:

\bigskip

\noindent {\bf Proposition~5.3.} {\sl Let $k$ be a field, 
let $1 \to N \to \Gamma \ {\buildrel f \over \to} \ G \to 1$ be a split short exact sequence of finite groups {\new with $N \ne 1$}, and let $Y^* \to X^* = \Spec k[[x,t]]$ be a $G$-Galois connected normal branched cover.  

a) Then there is a $\Gamma$-Galois connected normal branched cover $Z^* \to X^*$ dominating $Y^* \to X^*$.

b) The cover $Z^* \to X^*$ may be chosen such that $Z^* \to Y^*$ is totally split at the generic points of the ramification locus of $Y^* \to X^*$.

c) The set of isomorphism classes of such covers $Z^* \to X^*$ has cardinality equal to that of $k((x,t))$ (whether or not condition (b) is required){\new; and moreover there is a linearly disjoint set of such covers having the same cardinality}.}

\medskip

\noindent{\it Proof.} After a change of variables of the form $x'=x, t'=t+x^n$ for some $n>0$, we may assume that $Y^* \to X^*$ is unramified over the generic point of $(t=0)$.

By Lemma~5.2, we may extend $Y^* \to X^*$ to a cover $\bar Y \to \bar X$ of the $x$-line $\bar X = \PP^1_{k[[t]]}$ over $k[[t]]$.  
Let $\eta$ be the unique closed point of $Y^*$, and also write $\eta$ for the image of this point in $\bar Y$; this is the unique closed point of $\bar Y$ where $x=0$.
In order to prove the proposition, we will use Theorem~4.1 to construct an appropriate $\Gamma$-Galois cover of ${\bar{X}}$ dominating ${\bar{Y}}$, and then restrict that to obtain the desired $\Gamma$-Galois cover of $X^*$ dominating $Y^*$.  

Let $Y \to X$ be the general fibre of $\bar Y \to \bar X$, and let
$K$ be the function field of $X$ (or equivalently, of $\bar X$).  
Let $\alpha:G_K \to G$ be the surjection corresponding to the $G$-Galois cover $Y \to X$, and consider the finite split embedding problem ${\cal E}= (\alpha:G_K \to G, f:\Gamma \to G)$.  
To give a proper solution $\beta:G_K \to \Gamma$ to ${\cal E}$ is equivalent to giving a connected normal $\Gamma$-Galois cover $Z \to X$ dominating $Y \to X$; and for such a cover, we may consider the normalization $\bar Z$ of $\bar X$ in $Z$.  In that situation, there are then intermediate $N$-Galois covers $Z \to Y$ and $\bar{Z} \to \bar{Y}$ of the $\Gamma$-Galois covers $Z \to X$ and $\bar{Z} \to \bar{X}$.

Since $\bar X = \PP^1_{k[[t]]}$ is smooth over $k[[t]]$, we may apply Theorem~4.1, taking $D \subset Y$ to be the ramification locus of $Y \to X$.  By that result, we know in fact that there exists a proper solution $\beta:G_K \to \Gamma$ to $\cal E$ such that the intermediate cover $Z \to Y$ is totally split over the ramified points of $Y \to X$ (by part (b)), and such that $\bar Z$ is totally ramified at the unique point of $\bar Y$
over $\eta$ (by part (d) of Theorem~4.1; there $\eta=Q$).  

Let $Z^* \to X^*$ be the pullback of ${\bar{Z}} \to {\bar{X}}$ via $X^* \to {\bar{X}}$; this is a $\Gamma$-Galois cover.
Since ${\bar {Z}} \to {\bar{X}} $ dominates ${\bar{Y}} \to {\bar{X}}$, we see 
that $Z^* \to X^*$ dominates $Y^* \to X^*$.
Since $\bar{Z} \to \bar{Y}$ is totally ramified over $\eta$, the same is true for $Z^* \to Y^*$.  But $\eta$ is the unique closed point of $Y^*$.  So $Z^*$ has only one closed point, and is therefore connected.  Being normal, it is irreducible.  Here $Z^* \to Y^*$ splits completely over the generic points of the ramification locus of $Y^* \to X^*$, because this is a local condition and this
splitting holds on the generic fibre $Z \to Y$ of
${\bar Z} \to {\bar Y}$.  So $Z^* \to X^*$ is as asserted.

By Theorem~4.1, parts (c) and (d), there were ${\rm card}\, k((t))$ non-isomorphic choices above for $Z \to X$, and there were the same number of pullbacks $Z^* \to X^*$.  
{\new Moreover that theorem asserted there is a subset of choices for $Z \to X$ such that the intermediate covers $Z \to Y$ are linearly disjoint, as are the pullbacks $Z^* \to Y^*$; and that the sets of these covers and these pullbacks each have cardinality ${\rm card}\, k((t))$.}
But 
${\rm card}\, k((t)) = {\rm card}\, k[[t]] = ({\rm card}\, k)^{\aleph_0} = {\rm card}\, k[[x,t]] = {\rm card}\, k((x,t))$ (where, as usual, $\aleph_0$ denotes ${\rm card}\, \ZZ$).  So there is the asserted cardinality for the covers $Z^* \to X^*$ satisfying (a) and (b) of the proposition{\new, and the asserted {\new independence} condition holds}.  The same is true for covers $Z^* \to X^*$ that are just required to satisfy (a) of the proposition, by applying Theorem~4.1 with $D = \emptyset$ there, instead of taking $D$ to be the ramification locus of $Y \to X$ as before.
\qed

\bigskip

Restating Proposition~5.3, we obtain our main result, Theorem~5.1:

\bigskip

\noindent{\it Proof of 5.1.}  Given any finite split embedding problem 
${\cal E}= (\alpha:G_K \to G, f:\Gamma \to G)$ for $G_K$, consider the $G$-Galois connected normal cover $Y^* \to X^* = \Spec k[[x,t]]$ corresponding to $f$.  By Proposition~5.3, there are precisely ${\rm card}\, k[[x,t]] = {\rm card}\, K$ non-isomorphic $\Gamma$-Galois covers $Z^* \to X^*$ that dominate $Y^* \to X^*$.  These covers correspond in turn to ${\rm card}\, K$ distinct proper solutions to the embedding problem ${\cal E}$.  
{\new Moreover there is a linearly disjoint subset of these covers having the same cardinality, corresponding to ${\rm card}\, K$ independent proper solutions to ${\cal E}$.  So the absolute Galois group of $K$ is {\new semi-free} of rank equal to ${\rm card}\, K$.}
\qed

\

\bigskip

\noindent{\bf References.}

\medskip

\noindent [AGV] M.~Artin, A.~Grothendieck, J.-L.~Verdier.  ``Th\'eorie des topos et cohomologie \'etale des sch\'emas'' (SGA 4, vol.\ 3). Lecture Notes in Mathematics {\bf 305},
Springer-Verlag, 1973.

\smallskip

\noindent [DD] P.~D\`ebes, B.~Deschamps.  The regular inverse Galois problem over large fields.  In ``Geometric Galois actions'', vol.~2 (L.~Schneps and P.~Lochak, eds.), London Math.\ Soc.\ Lec.\ Note Ser.\ {\bf 243}, Cambridge U.\ Press, 1997, pp.~119-138.

\smallskip

\noindent [FJ] M.~Fried, M.~Jarden. ``Field Arithmetic.'' Ergebnisse Math.\ series, {\bf 11}, Springer-Verlag, 1986.

\smallskip

\noindent [GMP] B.~Green, M.~Matignon, F.~Pop.  On valued function fields II: Regular functions and elements with the uniqueness property.  J.\ Reine angew.\ Math., {\bf 412} (1990), 128-149.

\smallskip

\noindent [Gr61] A.~Grothendieck. ``\'El\'ements de G\'eom\'etrie Alg\'ebrique'' (EGA) III, $1^e$ partie, Publ.\ Math.\ IHES, vol.~11 (1961).

\smallskip

\noindent [Gr66] A.~Grothendieck. ``\'El\'ements de G\'eom\'etrie Alg\'ebrique'' (EGA) IV, $3^e$ partie, Publ.\ Math.\ IHES, vol.~28 (1966).

\smallskip

\noindent [Gru] K.W.~Gruenberg.   Projective profinite groups.  J.\ London Math.\ Soc., {\bf 42} (1967), 155-165.

\smallskip

\noindent[HJ] D.~Haran, M.~Jarden.  Regular split embedding problems over complete valued fields.  Forum Mathematicum, {\bf 10} (1998), 329-351.

\smallskip

\noindent [Ha87] D.~Harbater.  Galois coverings of the arithmetic line.  In ``Number Theory:  New York, 1984-85'' (D.V.\ and G.V.~Chudnovksy, eds.), Springer LNM {\bf 1240} (1987), pp. 165-195.

\smallskip

\noindent[Ha95] D.~Harbater.  Fundamental groups and embedding problems in characteristic $p$.  In    
``Recent Developments in the Inverse Galois Problem'' (M.~Fried, et al., eds.), AMS Contemporary Mathematics Series {\bf 186}, 1995, pp.353-369.

\smallskip

\noindent[Ha02] D.~Harbater. Shafarevich conjecture.  In Supplement III, Encyclopaedia of Mathematics.  Managing Editor: M. Hazewinkel, Kluwer Academic Publishers, 2002, pp.360-361.

\smallskip

\noindent [Ha03] D.~Harbater.  Patching and Galois theory.  In ``Galois Groups and Fundamental Groups" (L. Schneps, ed.), MSRI Publications series, {\bf 41}, Cambridge University Press, 2003, pp.313-424.

\smallskip

{\new 

\noindent [Ha08] D.~Harbater.  On function fields with free absolute Galois groups. {\new To appear in Crelle's Journal (online reference DOI 10.1515/CRELLE.2009.051).}  Also available at 
 
\noindent{\tt http://www.math.upenn.edu/{$\sim$}harbater/freegal.pdf}

\smallskip}

\noindent [HS99] D.~Harbater, K.~Stevenson. Patching and Thickening Problems.  J.\ Algebra, {\bf 212} (1999), 272-304.

\smallskip
 
\noindent[HS03] D.~Harbater, K.~Stevenson.  Abhyankar's local conjecture on fundamental groups.  In ``Algebra, Arithmetic and Geometry with Applications'' (Abhyankar 70th 
birthday conference proceedings, C.Christensen et al., eds.), 
Springer-Verlag, 2003, pp. 473-485.

{\smallskip

\new \noindent[HS05] D.~Harbater, K.~Stevenson,
Local Galois theory in dimension two. Advances in 
Math.\ (special issue in honor of M.~Artin's 70th birthday), 
{\bf 198} (2005), 623-653.

\smallskip

\noindent[HS06] D.~Harbater, K.~Stevenson, Local Galois theory in dimension two, Oberwolfach Report No.~6 (2006), 341-344.}

\smallskip

\noindent[HP] D.~Harbater, M.~van der Put, with an appendix by R.~Guralnick.  Valued fields and covers in characteristic $p$.  In ``Valuation Theory and its Applications,'' Fields Institute Communications, vol.~32, ed.\ by F.-V.~Kuhlmann, S.~Kuhlmann and M.~Marshall, 2002, pp.175-204.

\smallskip

\noindent [Hrt] R. Hartshorne. ``Algebraic geometry.'' Graduate Texts in Mathematics, vol.~52. Springer-Verlag, 1977.

\smallskip

\noindent [Iw] K.~Iwasawa.  On solvable extensions of algebraic number fields.  Annals of Math., {\bf 58} (1953), 548-572.

\smallskip

\noindent [Ja] M.~Jarden. On free profinite groups of uncountable rank. 
In ``Recent developments in the inverse
Galois problem'' (M.~Fried, ed.), AMS Contemporary Mathematics Series, {\bf 186}, 1995, pp.~371-383.

\smallskip

\noindent[Ka] N.~Katz.  Local-to-global extensions of representations of fundamental groups.
  Ann.\ l'inst.\ Fourier, {\bf 36} (1986), 69-106.

\smallskip

\noindent [La]   S.~Lang.  ``Algebraic Number Theory''.  Addison-Wesley, 1970.

\smallskip

\noindent [Le] T.~Lefcourt.   Galois groups and complete domains.  Israel J.\ Math.\ {\bf 114} (1999), 323-346.

\smallskip

\noindent[Po95] F.~Pop.  \'Etale Galois covers of affine smooth curves.  Invent.\ Math., {\bf 120} (1995), 555-578.

\smallskip

\noindent[Po96] F.~Pop.  Embedding problems over large fields.  Ann.\ Math., {\bf 144} (1996), 1-34.

\smallskip

\noindent [Pr] R.~Pries. Construction of covers with formal and rigid geometry.  In ``Courbes semi-stables et groupe fondamental en g\'eom\'etrie alg\'ebrique'' (J.-B.~Bost, F.~Loeser, M.~Raynaud, eds.), Progress in Math.\ {\bf 187}, Birkh\"auser, 2000, pp.~157-167.

\smallskip

\noindent[Se] J.-P.~Serre. ``Cohomologie Galoisienne''.
Lecture Notes in Mathematics, {\bf 5}, Springer-Verlag, 1964.

\bigskip

\small

\noindent Author information: 

\medskip

\noindent David Harbater: Department of Mathematics,
University of Pennsylvania, Philadelphia, PA 19104-6395.
{\smit E-mail address}: {\smtt harbater@math.upenn.edu}

\medskip

\noindent Katherine F.\ Stevenson: Dept.\ of Mathematics,
California State University at Northridge, Northridge, CA 91330.
{\smit E-mail address}: {\smtt katherine.stevenson@csun.edu}

\end